\documentclass[3p,11pt,authoryear]{elsarticle}

\usepackage[utf8]{inputenc}
\usepackage{graphicx}
\usepackage{hyperref}    
\usepackage{multirow}

\usepackage{marginnote}
\usepackage{amstext,psfrag}
\usepackage[dvips]{epsfig}
\usepackage{ragged2e}  
\usepackage{gensymb}
\usepackage{graphicx}
\usepackage{epstopdf}
\usepackage{color}
\usepackage{latexsym}
\usepackage{float}
\usepackage{array}
\usepackage{tabu}
\usepackage{multirow}
\usepackage{textcomp}
\usepackage{float}
\usepackage{subfig}
\usepackage{amsmath}
\usepackage{amsthm}
\usepackage{amssymb}
\usepackage[version=4]{mhchem}
\usepackage{siunitx}
\usepackage{chngcntr}
\usepackage{apptools}
\usepackage{longtable,tabularx}
\usepackage{newcmd}

\newtheorem{Definition}{Definition}
\newtheorem{Theorem}{Theorem}
\newtheorem{Corollary}{Corollary}
\newtheorem{Proposition}{Proposition}
\newtheorem{Lemma}{Lemma}
\newtheorem*{Proof*}{Proof}

\newtheorem{Remark}{Remark}

\begin{document}

\begin{frontmatter}
	
	
	
	\title{Geometric Extended State Observer on \SE \ with Fast Finite-Time Stability: Theory and Validation on a Rotorcraft Aerial Vehicle}
	
	
	\author{Ningshan Wang}
	\author{Reza Hamrah}
	\author{Amit K. Sanyal}
	\author{Mark N. Glauser}
	
	\affiliation{organization={Department of Mechanical \& Aerospace Engineering, Syracuse University},
		city={Syracuse},
		state={NY 13244},
		country={US}}
	
	\begin{abstract}
		This article presents an extended state observer for 
		vehicle modeled as a rigid body in three-dimensional translational and rotational motions. The extended state observer is applicable to a rotorcraft aerial vehicle with a fixed plane of rotors, modeled as an under-actuated system on the tangent bundle of the six-dimensional Lie group of rigid body motions, $\SE$. The extended state observer is designed to estimate the resultant external disturbance force and disturbance torque acting on the vehicle. It guarantees stable convergence of 
		disturbance estimation errors in finite time when the disturbances are constant and finite time convergence to a bounded neighborhood of zero errors for time-varying disturbances. This extended state observer design is based on a H\"{o}lder-continuous fast finite time stable differentiator that is similar to the super-twisting algorithm, to obtain fast 
		convergence. Numerical simulations are conducted to validate the proposed extended state observer. The proposed extended state observer is compared with other existing research to show its advantages. A set of experimental results implementing disturbance rejection control using feedback of disturbance estimates from the extended state observer is also presented.
	\end{abstract}
	
	
	
	
\end{frontmatter}

\section{Introduction}\label{sec:Intro}

Rotorcraft unmanned aerial vehicles (UAVs) are increasingly being used in various applications, such as security and monitoring, infrastructure inspection, agriculture, wildland management, package delivery, and remote sensing. However, these UAVs are frequently exposed to dynamic uncertainties and disturbances caused by turbulence induced by airflow around structures or regions. Therefore, it is crucial to ensure robust flight control performance in such challenging environments, with guaranteed stability margins even in the presence of dynamic disturbances and uncertainties. 

Recent research articles on rotorcraft UAV tracking control have used various methods to tackle the adverse effects of disturbances and uncertainties during the flight. \cite{torrente2021data} use  Gaussian processes to complement the nominal dynamics of the multi-rotor in a model predictive control (MPC) pipeline. \cite{hanover2021performance} use an explicit scheme to discretize the dynamics for a nonlinear MPC solved by optimization. 
\cite{bangura2017thrust} use the propeller aerodynamics as a direct feedforward term on the desired thrust to re-regulate the thrust command of the rotors. \cite{craig2020geometric} implement a set of pitot tubes onto the multi-rotor aircraft to directly sense the aircraft's airspeed. 
With the knowledge of propeller aerodynamic characteristics, the airspeed is then utilized to obtain the disturbance forces and torques as feedforward terms to enhance control performance. \cite{bisheban2020geometric} use artificial neural networks to obtain disturbance forces and torques with the kinematics information of the aircraft and then use the baseline control scheme based on the article by \cite{lee2010geometric} in their tracking control scheme design. The methods used in these research articles either need high computational efforts \cite{torrente2021data, hanover2021performance,bisheban2020geometric} or require precise modeling of the aerodynamic characteristics of the rotorcraft propellers \cite{bangura2017thrust,craig2020geometric}, to obtain satisfactory control performance against disturbances.

In this article, we develop and use extended state observers to  estimate disturbance force and disturbance torque vectors acting on a rotorcraft UAV. Extended state observers, 
along with disturbance observers and unknown input observers, are commonly used in association with the robust control technique known as {\em active disturbance rejection control} (ADRC), which can be traced back to the dissertation by \cite{hartlieb1956cancellation}.
In an ADRC scheme, estimates of unknown disturbance inputs from a disturbance observer (DO) or an extended state observer (ESO) are first obtained and then utilized in the control design to reject the disturbance. 
See, for example, \cite{huang2001flight, shao2018robust, mechali2021observer, cui2021adaptive} on application of ESO and \cite{chen2003nonlinear, liu2022fixed, bhale2022finite, sanyal2022discrete} on application of DO in ADRC schemes. \cite{jia2022accurate} employ the disturbance model obtained by \cite{faessler2017differential} and then estimate the drag coefficient as a parameter. This disturbance model is also employed by \cite{moeini2021exponentially}. 

There are several methods to ensure the stability of ESO/DO designs used for rotorcraft tracking control. The linear ESO by \cite{shao2018robust} is asymptotically stable (AS). \cite{mechali2021observer} use the concept of geometric homogeneity \cite{rosier1992homogeneous} to obtain an FTS ESO. A similar method is proposed in the ESO design by \cite{guo2011convergence}. The 
Lyapunov functions/candidates used in the ESO stability analysis by \cite{mechali2021observer} and \cite{guo2011convergence} are based on \cite{rosier1992homogeneous}, and are implicit. \cite{jia2022accurate}, \cite{moeini2021exponentially} and \cite{liu2022fixed} use variants of the DO 
proposed by \cite{chen2003nonlinear}. Another approach is to use the super-twisting algorithm (STA) \cite{moreno2012strict} to design ESO. \cite{xia2010attitude} use this method in ESO design for spacecraft attitude control, and \cite{cui2021adaptive} design an adaptive super-twisting ESO using a similar method for an ADRC scheme applied to rotorcraft UAV.

In much of the prior literature for rotorcraft UAV attitude control with ESO/DO for disturbance torque estimation and rejection in rotational dynamics, the attitude kinematics of the ESO/DO are either based on local 
linearization or represented using local coordinates (like Euler angles) or quaternions. Local coordinate representations can have singularity issues  (e.g., gimbal lock with Euler angles), while quaternion 
representations may cause instability due to unwinding \cite{bhat2000topological,chaturvedi2011rigid}. In situations where the UAVs have to carry out aggressive maneuvers, as in rapid collision avoidance for example, disturbance estimation and rejection from such schemes may not be reliable or accurate enough for precise control of the UAV.

This article presents an ESO on $\SE$ for rotorcraft UAVs to provide reliable disturbance estimation under complex and challenging aerodynamic environments. The ESO on \SE \ estimates the disturbance forces and torques during the flight of a UAV in both translational and rotational dimension. The proposed ESO is fast finite-time stable (FFTS), abbreviated as FFTS-ESO. This FFTS-ESO design is based on a novel H\"{o}lder-continuous fast finite-time stable differentiator (HC-FFTSD). We carry out several sets of numerical simulations to show the validity of the proposed FFTS-ESO. Moreover, the proposed FFTS-ESO is compared with other existing research to show its advantages in the conducted simulations. A set of experimental results implementing a disturbance rejection mechanism using feedback of disturbance estimates from the FFTS-ESO is also presented. In the experiment, we hover the UAV in front of the turbulent flows generated by a fan array wind tunnel (FAWT). We obtain statistical information from the hot-wire measurements on the turbulent incoming flows. We observe the pose of the UAV to evaluate its flight control performance.

We highlight some unique contributions of this article.
\begin{itemize}
	\item The proposed ESO is the major contribution of this article. The pose of the rotorcraft is represented directly on the Lie group of rigid body transformations, the special Euclidean group $\SE$. Unlike the ESO and DO designs reported by \cite{mechali2021observer}, \cite{shao2018robust}, and \cite{cui2021adaptive}, which use Euler angles or quaternions for attitude representation or do not include attitude kinematics, like the DO by \cite{bhale2022finite} in disturbance torque estimation, the pose of the aircraft in this article is represented in \SE to avoid kinematic singularities. We do not use local coordinates (like Euler angles) or (dual) quaternions for pose representation so that we avoid singularities due to local coordinate representations or quaternion unwinding, as reported by \cite{bhat2000topological}, and \cite{chaturvedi2011rigid}.  To the best of the author's knowledge, there is no existing publication on aircraft disturbance observation using ESO with pose representation on $\SE$. 
	\item The proposed FFTS-ESO is based on the HC-FFTSD. The commonly used geometric homogeneity method \cite{rosier1992homogeneous,guo2011convergence,liu2019state,wang2021holder,wang2022holder}, cannot provide a straightforward (or explicit) Lyapunov function to prove the finite-time stability of the scheme. The (implicit) form of their Lyapunov functions is by \cite{rosier1992homogeneous}. This implicit Lyapunov function complicates the robustness analysis under measurement noise and time-varying disturbances when that analysis is essential for an ESO designed for disturbance estimation in ADRC schemes. We propose HC-FFTSD as an approach inspired by the STA \cite{moreno2012strict,vidal2016output} of sliding-mode control (SMC). This approach gives a straightforward design of a strict Lyapunov function, which is explicit, and therefore avoids the weakness mentioned above.
	\item Based on the HC-FFTSD, the proposed FFTS-ESO schemes are both FFTS and H\textup{\"{o}}lder-continuous, unlike the common STA and other FTS schemes that use discontinuous methods like terminal sliding-mode. Therefore, the proposed FFTS-ESO avoids the potentially harmful chattering phenomenon \cite{sanyal2015finite}, while maintaining FTS convergence. 
	\item With explicit Lyapunov function in the stability analysis, we present proof of the robustness of the proposed FFTS-ESO under time-varying disturbing forces, torques, and measurement noise. To the best of the authors' knowledge, there is no prior research on the noise robustness of ESO using Lyapunov analysis.
\end{itemize}
The remainder of the article is as follows. Section \ref{sec:Preliminary} presents some preliminary results that are needed to obtain sufficient conditions for the stability of the ESO and ADRC schemes. HC-FFTSD is presented, along with its stability and robustness analysis in Section \ref{sec:Differentiator}. In Section \ref{sec:Problem}, the ESO design problem is formulated. Section \ref{sec:ESO} describes the detailed FFTS-ESO design, which is based on the differentiator design in Section \ref{sec:Differentiator}. Numerical simulations are conducted in Section \ref{sec:Numerical}. Section \ref{sec:Experiment} describes the conducted UAV flight experiment with the UAV exposed to the disturbances generated by the FAWT in details. We conclude the paper, in Section \ref{sec:Conclusion}, by summarizing the results and highlighting directions for forthcoming research.

\section{Preliminaries}\label{sec:Preliminary}
The statements and definitions in this section are used in the technical results obtained in later sections. The statements given here give the conditions under which a continuous time system is finite-time stable, fast finite-time stable, and practically finite-time stable using Lyapunov analysis, and the last statement is used in developing the main result. 
\begin{Lemma}[Finite-time stable]\label{lem:FTS}
	\textup{\cite{bhat2000finite}} Consider the following system of differential equations,
	\begin{align}\label{eqn:Nonlinear System}
		\dot{x}(t) = f(x(t)),\ f(0)=0,\ x(0)=x_0,
	\end{align}
	where $f: \cD \rightarrow\bR^n$ is continuous on an open neighborhood $\cD \subset \bR^n$ of the origin, and let there be a continuous and differentiable function $V(x(t))$ that is positive definite. Let the time derivative of $V(x)$ satisfy the following inequality:
	\begin{align}\label{eqn:Preliminary FTS Lyapunov Inequality}
		\dot{V}\leq -\lambda V^{\alpha},
	\end{align}
	where $x(t)\in \cD \backslash \{0\}$, $\lambda>0$, $\alpha\in]0,1[$. Then the system \eqref{eqn:Nonlinear System} is FTS at the origin, which means $\forall x_0\in \cD$, $x$ can reach the origin in finite time. Moreover, the settling time $T$, the time needed to reach the origin, satisfies
	\begin{align}\label{eqn:Preliminary FTS Settling Time}
		T \leq \frac{V^{1-\alpha}(x_0)}{\lambda(1-\alpha)}.
	\end{align}
\end{Lemma}

\begin{Lemma}[Fast finite-time stable]\label{lem:FFTS}
	\textup{\cite{yu2005continuous}} Consider the system \eqref{eqn:Nonlinear System} and let there be a continuous and differentiable function $V(x(t))$ that is positive definite. Let the time derivative of $V(x)$ satisfy the following inequality:
	\begin{align}\label{eqn:Preliminary FFTS Lyapunov Inequality}
		\dot{V} \leq -\lambda_1 V-\lambda_2 V^{\alpha},
	\end{align}
	where $x(t)\in \cD \backslash \{0\}$, $\lambda_1, \lambda_2>0$, $\alpha\in]0,1[$. Then the system \eqref{eqn:Nonlinear System} is FFTS at the origin and the settling time $T$ satisfies:
	\begin{align}\label{eqn:Preliminary FFTS Settling Time}
		T \leq \frac{1}{\lambda_1(1-\alpha)}\textup{ln}\frac{\lambda_1V^{1-\alpha}(x_0)+\lambda_2}{\lambda_2}.
	\end{align}
\end{Lemma}

\begin{Lemma}[Practically finite-time stable]\label{lem:PFTS}
	\textup{\cite{yu2005continuous,zhu2011attitude}} Consider the system \eqref{eqn:Nonlinear System} and let there be a continuous and differentiable function $V(x)$ that is positive definite. Let the time derivative of $V(x)$ satisfy the following inequality:
	\begin{align}\label{eqn:Preliminary PFTS Lyapunov Inequality}
		\dot{V} \leq -\lambda_1 V-\lambda_2 V^{\alpha} +\eta,
	\end{align}
	with $x(t)\in \cD \backslash \{0\}$,  $\lambda_1, \lambda_2>0$, and $\alpha\in]0,1[$. Then the system \eqref{eqn:Nonlinear System} is practical finite-time stable (PFTS) at the origin, which means that
	the solution of \eqref{eqn:Nonlinear System} will converge to the following set in finite time
	\begin{align*}
		\left\{x\,  \bigg| \, V(x)\leq \textup{min}\left\{ \frac{\eta}{(1-\theta_0)\lambda_1}, \left(\frac{\eta}{(1-\theta_0)\lambda_2} \right)^{\frac{1}{\alpha}} \right\} \right\},
	\end{align*}
	where $0<\theta_0<1$. The settling time $T$ is bounded above as follows:
	\begin{align*}
		T \leq  \textup{max}&\left\{ t_0  + \frac{1}{\theta_0\lambda_1(1-\alpha)}\textup{ln}\frac{\theta_0\lambda_1V^{1-\alpha}(x_0)+\lambda_2}{\lambda_2},  t_0  + \frac{1}{\lambda_1(1-\alpha)}\textup{ln}\frac{\lambda_1V^{1-\alpha}(x_0)+\theta_0\lambda_2}{\theta_0\lambda_2}\right\}.
	\end{align*}
\end{Lemma}

\begin{Lemma}\label{lem:Binom}
	\textup{\cite{hardy1952inequalities}}Let $x$ and $y$ be non-negative real numbers and let $p\in ]1,2[$. Then 
	\begin{align}\label{eqn:Bires}
		x^{\frac{1}{p}}+ y^{\frac{1}{p}} \ge (x+y)^{\frac{1}{p}}.
	\end{align}
	Moreover, the above inequality is a strict inequality if both $x$ and $y$ are non-zero. 
\end{Lemma}

\begin{Definition}\label{def:H}
	Define $H: \mathbb{R}^3 \times \mathbb{R}\rightarrow \textup{Sym(3)}$, the space of symmetric $3\times3$ matrices, as follows:
	\begin{equation}\label{eqn:H}
		H(x,k) := I - \frac{2k}{x\Tp x} xx\Tp. 
	\end{equation}
	
\end{Definition}

\begin{Lemma}\label{lem:Inequality Noise Robustness}
	Let $\mu\in\bR^n/\{0\}$ and $\alpha\in]0,1/2[$. Consider $\cD:\bR^n\setminus \{0,-\mu\}$ and define $\phi(x):\cD \rightarrow \bR^+$ as:
	\begin{align}\label{eqn:Function phi}
		\begin{split}
			&\phi(x) := Y(x)\Tp Y(x), \mbox{ where } \\
			&Y(x) :={\|x\|^{-2\alpha}}x - {\|x+\mu\|^{-2\alpha}}(x+\mu).
		\end{split}
	\end{align}
	The global maximum of $\phi(x)$ is at $x=-\mu/2$.
\end{Lemma}
We provide the proof of Lemma \ref{lem:Inequality Noise Robustness} in the appendix. 

\section{H\"{o}lder-Continuous Fast Finite-Time Stable Differentiator (HC-FFTSD)}\label{sec:Differentiator}
In this section, we design the error dynamics for the proposed ESO in Section \ref{sec:ESO} in the form of an HC-FFTSD. We analyze the stability and robustness of the proposed HC-FFTSD in this section, to support the development of the ESO design in Section \ref{sec:ESO}. Theorem \ref{thr:FFTS Differentiator} gives the proposed HC-FFTSD with its stability properties. Corollary \ref{cor:FTS Differentiator Disturbance Robustness} describes the convergence performance of the differentiator under external disturbances.  Corollary \ref{cor:FTS Differentiator Noise Robustness} describes the convergence performance of the differentiator under measurement noise. In the analysis that follows, $e_1 \in \bR^n$ stands for the measurement estimation error and $e_2 \in \bR^n$ stands for the disturbance estimation error in the ESO error dynamics, respectively. In this section and the remainder of this paper, we denote the minimum and maximum eigenvalues of a matrix by $\lambda_{\min}(\cdot)$ and $\lambda_{\max}(\cdot)$, respectively.
\begin{Theorem}\label{thr:FFTS Differentiator}
	Let $p\in]1,2[$ and $k_3>0$. Define $\phi_1(\cdot): \bR^n \rightarrow \bR^n$ and  $\phi_2(\cdot): \bR^n \rightarrow \bR^n$ as follows: 
	\begin{align}\label{eqn:phi1phi2}
		\begin{split}
			\phi_1(e_1) &=    k_3 e_1 +(e_1\Tp e_1)^\frac{1-p}{3p-2}e_1, \\
			\phi_2(e_1) &=    k_3^2 e_1 + \frac{2k_3(2p-1)}{3p-2} (e_1\Tp e_1)^\frac{1-p}{3p-2}e_1+\frac{p}{3p-2}(e_1\Tp e_1)^\frac{2(1-p)}{3p-2}e_1. 
		\end{split}
	\end{align}
	Define the differentiator gains $k_1, k_2>0$ and  $\mathcal{A}^*\in \bR^{2\times2}$, as:
	\begin{align}\label{eqn:A}
		\mathcal{A}^* = 
		\left[
		\begin{array}{cc}
			-k_1  &  1 \\
			-k_2  & 0 \\
		\end{array}
		\right],	
	\end{align}
	which makes $\mathcal{A}^*$ a Hurwitz matrix. Thereafter, the differentiator design:
	\begin{align}\label{eqn:Differentiator}
		\begin{array}{ll}
			\dot{e}_1 &= -k_1 \phi_1(e_1) + e_2, \\ 
			\dot{e}_2 &= -k_2 \phi_2(e_1), 
		\end{array}
	\end{align}
	ensures that $(e_1\Tp,e_2\Tp)\in\bR^{2n}$ converges to the origin in a fast finite-time stable manner.
\end{Theorem}
\begin{Proof*}
	{\em
		The proof of Theorem \ref{thr:FFTS Differentiator} is based on Theorem 1 by \cite{vidal2016output}, and Theorem 1 by \cite{moreno2012strict}. Two properties of $\phi_1$ and $ \phi_2$ are provided as follows. \\
		\textit{Property 1 (P1): The Jacobian of $\phi_1(e_1)$, denoted $\phi_1'(e_1)$, is given as follows:}
		\begin{align}\label{eqn:Property1A}
			\begin{split}
				\phi'_1(e_1)&=\frac{\textup{d}\phi_1(e_1)}{\textup{d}e_1} = k_3 I +(e_1\Tp e_1)^\frac{1-p}{3p-2}
				\bigg[
				I-\frac{2(p-1)}{3p-2}\frac{e_1e_1\Tp}{e_1\Tp e_1}
				\bigg],
			\end{split}
		\end{align}
		\textit{so that the following identity holds:}
		\begin{align}\label{eqn:Property1B}
			\phi_2(e_1)=\phi'_1(e_1)\phi_1(e_1)
		\end{align}
		\textit{Property 2 (P2): $\phi_1'$ is a positive definite matrix, which means $\forall w\in \bR^{2n}, e_1\in \bR^n$,}
		\begin{align}\label{eqn:Property2A}
			\lambda_{\textup{min}} \{\phi’_1(e_1)\}||w||^2  \leq w\Tp\phi'_1(e_1)w\leq \lambda_{\textup{max}} \{\phi’_1(e_1)\}||w||^2. 
		\end{align}
		\textit{The maximum and minimum eigenvalues of $\phi'_1(e_1)$ employed in \eqref{eqn:Property2A} are as given below: }
		\begin{align}
			&\lambda_{\textup{max}} \{\phi'_1(e_1)\}=k_3 + (e_1\Tp e_1)^{\frac{1-p}{3p-2}}, \label{eqn:Property2B max} \\
			&\lambda_{\textup{min}} \{\phi'_1(e_1)\}=k_3 + (e_1\Tp e_1)^{\frac{1-p}{3p-2}} \frac{p}{3p-2}\label{eqn:Property2B min}. 
		\end{align}
		From Theorem 5.5 by Chen~\cite{chen1984linear}, we know that for a Hurwitz matrix $\mathcal{A}^*$ as in \eqref{eqn:A}, $\forall\, \mathcal{Q}^* \in \bR^{ 2\times 2}$ where $\mathcal{Q}^*  \succ 0$, the Lyapunov equation: 
		\begin{align}\label{eqn:Lyapunov Equation}
			(\mathcal{A}^*)\T \mathcal{P}^*+\mathcal{P}^*\mathcal{A}^* = -\mathcal{Q}^* ,
		\end{align}
		has a unique solution $\mathcal{P}^*\succ 0$. Express the positive definite matrices $\mathcal{P}^*$ and $\mathcal{Q}^*$ in components as:
		\begin{align}
			\mathcal{P}^* = 
			\left[
			\begin{array}{cc}
				p_{11}  &  p_{12} \\
				p_{12}  &  p_{22} \\
			\end{array}
			\right], \
			\mathcal{Q}^* = 
			\left[
			\begin{array}{cc}
				q_{11}  &  q_{12} \\
				q_{12}  &  q_{22} \\
			\end{array}
			\right]. \notag
		\end{align}
		With $\mathcal{P}^*$ defined as the solution to \eqref{eqn:Lyapunov Equation},  $\mathcal{A}^*$, 
		$\mathcal{P}^*$ and $\mathcal{Q}^*$ can be augmented to $\mathcal{A}, \mathcal{P}, \mathcal{Q} \in \bR^{2n\times 2n}$, as follows:
		\begin{align*}
			\begin{split}
				&\mathcal{A} = 
				\left[
				\begin{array}{cc}
					-k_1 I  &  I \\
					-k_2 I &  0 \\
				\end{array}
				\right], 
				\mathcal{P} = 
				\left[
				\begin{array}{cc}
					p_{11}I  &  p_{12}I \\
					p_{12}I &  p_{22}I \\
				\end{array}
				\right], 
				\mathcal{Q} = 
				\left[
				\begin{array}{cc}
					q_{11}I  &  q_{12}I \\
					q_{12}I  &  q_{22}I \\
				\end{array}
				\right]. 
			\end{split}
		\end{align*}
		The augmented matrices $\mathcal{A}, \mathcal{P}, \mathcal{Q}$ defined above also satisfy a Lyapunov equation as given below:
		\begin{align}\label{eqn:Augmented Lyapunov Equation}
			\mathcal{A}\T\mathcal{P}+\mathcal{P}\mathcal{A}=-\mathcal{Q}.
		\end{align}
		Further, the eigenvalues of $\mathcal{P}$ and $\mathcal{P}^*$, are related such that $\lambda_{\textup{min}}\{\mathcal{P}^*\} = \lambda_{\textup{min}}\{\mathcal{P}\} $, and $\lambda_{\textup{max}}\{\mathcal{P}^*\} = \lambda_{\textup{max}}\{\mathcal{P}\} $. Similar relations hold for $\mathcal{Q}$ and $\mathcal{Q}^*$.
		Thus, with $\mathcal{P}$ as the solution to \eqref{eqn:Augmented Lyapunov Equation}, we consider the following Lyapunov candidate:
		\begin{align}\label{eqn:Lyapunov FTS Differentiator}
			V(e_1,e_2)=\zeta\T \mathcal{P}\zeta, 
		\end{align}
		where $\zeta \in \bR^{2n}$ is defined as $\zeta := [\phi_1\T(e_1),e_2\T]\T$ and $\mathcal{P}$ is the augmented $\mathcal{P}^*$, which is the unique solution of \eqref{eqn:Lyapunov Equation} for a given $\mathcal{Q}^*  \succ 0$. The upper and lower bounds of the Lyapunov candidate $V$ in \eqref{eqn:Lyapunov FTS Differentiator} are as given below:
		\begin{align}\label{eqn:Lyapunov Bound 1}
			\lambda_\textup{min}\left\{\mathcal{P}\right\} \|\zeta\|^2 \leq V(e_1,e_2) \leq \lambda_\textup{max}\left\{\mathcal{P}\right\} \|\zeta\|^2.
		\end{align}
		From \eqref{eqn:Lyapunov Bound 1}, we obtain the following two inequalities:
		\begin{align}
			&\lambda_\textup{min}\left\{\mathcal{P}\right\} (e\T_1e_1)^\frac{p}{3p-2} \leq  \lambda_\textup{min}\left\{\mathcal{P}\right\} \|\zeta\|^2 \leq V(e_1,e_2), \label{eqn:Lyapunov Bound 2} \\
			&k^2_3\lambda_\textup{min}\left\{\mathcal{P}\right\} e\T_1e_1 \leq  \lambda_\textup{min}\left\{\mathcal{P}\right\}\|\zeta\|^2 \leq V(e_1,e_2). \label{eqn:Lyapunov Bound 3}
		\end{align}
		$V(e_1,e_2)$ is differentiable everywhere except the subspace $\mathcal{S}=\{(e_1,e_2)\in\bR^{2n}|e_1 = 0 \}$.
		From \eqref{eqn:Differentiator} and Property (P1), we obtain the time derivative of $\zeta$ as follows,
		\begin{align}\label{eqn:zeta Derivtive}
			\begin{split}
				\dot{\zeta} &=
				\begin{bmatrix}
					\phi'_1(e_1)\dot{e}_1  \\
					\dot{e}_2 
				\end{bmatrix} =
				\begin{bmatrix}
					\phi'_1(e_1)(-k_1\phi_1(e_1)+e_2)\\
					-k_2  \phi'_1(e_1) \phi_1(e_1)
				\end{bmatrix} \\
				& = \mathcal{D}(e_1)\mathcal{A}\zeta ,
			\end{split}
		\end{align}
		where,
		\begin{align}\label{eqn:D}
			\begin{split}
				&\mathcal{D}(e_1) = \textup{diag}[\phi'_1(e_1),\phi'_1(e_1)]\in \bR^{2n\times 2n}, \\
				&\lambda_\textup{min}\left\{\mathcal{D}(e_1)\right\} = \lambda_\textup{min}\left\{\phi'_1(e_1)\right\}.
			\end{split}
		\end{align}
		With the expression of $\dot{\zeta}$ in \eqref{eqn:zeta Derivtive}, we obtain the time derivative of $V(e_1, e_2)$ as
		\begin{align}\label{eqn:Lyapunov Derivative 1}
			\begin{split}
				\dot{V} &= \dot{\zeta}\T \mathcal{P} \zeta + \zeta\T \mathcal{P}\dot{\zeta}  \\
				&= \zeta\T ((\mathcal{D}(e_1)\mathcal{A})\T \mathcal{P}+\mathcal{P}\mathcal{D}(e_1)\mathcal{A})\zeta  \\
				&= -\zeta\T \mathcal{\bar{Q}}(e_1) \zeta. 
			\end{split}
		\end{align}
		where $\mathcal{\bar{Q}}(e_1) $ is as
		\begin{align}\label{eqn:bar Q}
			\begin{split}
				\mathcal{\bar{Q}}(e_1) &= (\mathcal{D}(e_1)\mathcal{A})\T \mathcal{P}+\mathcal{P}\mathcal{D}(e_1)\mathcal{A} 
				=
				\begin{bmatrix}
					\mathcal{\bar{Q}}_{11}(e_1) & \mathcal{\bar{Q}}_{12}(e_1) \\
					\mathcal{\bar{Q}}_{12}(e_1) & \mathcal{\bar{Q}}_{22}(e_1) 
				\end{bmatrix}, \\
				\mathcal{\bar{Q}}_{11}(e_1) &= 2(k_1p_{11}+k_2p_{12})\phi'_1(e_1), \\
				\mathcal{\bar{Q}}_{12}(e_1) &= (k_1p_{12}+k_2p_{22} -p_{11})\phi'_1(e_1),\\
				\mathcal{\bar{Q}}_{22}(e_1) &= -2p_{12}\phi'_1(e_1). 
			\end{split}
		\end{align}
		With \eqref{eqn:bar Q} and \eqref{eqn:Augmented Lyapunov Equation}, we obtain $\mathcal{\bar{Q}} = \mathcal{Q}\mathcal{D}(e_1) $. Afterwards, with $\mathcal{Q},\; \mathcal{D}(e_1) \succ 0$, as defined in \eqref{eqn:Augmented Lyapunov Equation} and \eqref{eqn:D}, following inequality on their eigenvalues holds:
		With $\mathcal{Q} \succ 0$ and $\mathcal{D}(e_1)\succ 0$, we obtain following inequality on their eigenvalues,
		\begin{align}\label{eqn:Eigen Inequality 1}
			\lambda_\textup{min}\left\{\mathcal{Q} \mathcal{D}(e_1)\right\} \geq \lambda_\textup{min}\left\{\mathcal{Q} \right\} \lambda_\textup{min}\left\{\mathcal{D}(e_1)\right\} >0.
		\end{align}
		With Property 2, substituting \eqref{eqn:Eigen Inequality 1} into \eqref{eqn:Lyapunov Derivative 1}, we obtain
		\begin{align}\label{eqn:Lyapunov Derivative 2}
			\begin{split}
				\dot{V} &= -\zeta\T (\mathcal{Q} \mathcal{D}(e_1)) \zeta \\
				& \leq -\lambda_\textup{min}\left\{\mathcal{Q} \mathcal{D}(e_1)\right\}\zeta\T \zeta \\
				& \leq -\lambda_\textup{min}\left\{\mathcal{D}(e_1)\right\}\lambda_\textup{min}\left\{\mathcal{Q} \right\}\zeta\T \zeta 
			\end{split}
		\end{align}
		With $\lambda_\textup{min}\left\{\mathcal{D}(e_1)\right\} = \lambda_\textup{min}\left\{\phi'_1(e_1)\right\}$,  substituting \eqref{eqn:Property2B min} and \eqref{eqn:Lyapunov Bound 2} into \eqref{eqn:Lyapunov Derivative 2}, we obtain,
		\begin{align}\label{eqn:Lyapunov Derivative 3}
			\begin{split}
				\dot{V} & \leq-\Big[ k_3 + (e_1\T e_1)^{\frac{1-p}{3p-2}} \frac{p}{3p-2} \Big]\lambda_\textup{min}\left\{\mathcal{Q} \right\} \zeta\T \zeta \\
				&\leq -\frac{\lambda_\textup{min}\left\{\mathcal{Q} \right\}}{\lambda_\textup{max}\left\{\mathcal{P}\right\}}\Big[ k_3 + \left(\frac{V}{\lambda_\textup{min}\left\{\mathcal{P}\right\}}\right)^\frac{1-p}{p} \frac{p}{3p-2} \Big] V  \\
				& \leq -\gamma_1 V -\gamma_2 V^\frac{1}{p}, 
			\end{split}
		\end{align}
		where $\gamma_1$ and $\gamma_2$ are positive constants, defined as,
		\begin{align} \label{eqn:gamma1 gamma2}
			\begin{split}
				\gamma_1 &= k_3 \frac{\lambda_\textup{min}\left\{\mathcal{Q} \right\}}{\lambda_\textup{max}\left\{\mathcal{P}\right\}} = k_3 \frac{\lambda_\textup{min}\left\{\mathcal{Q}^* \right\}}{\lambda_\textup{max}\left\{\mathcal{P}^*\right\}};\\ 
				\gamma_2 &= \frac{\lambda_\textup{min}\left\{\mathcal{Q} \right\}\lambda_\textup{min}\left\{\mathcal{P}\right\}^\frac{p-1}{p}}{\lambda_\textup{max}\left\{\mathcal{P}\right\}} \frac{p}{3p-2} = \frac{\lambda_\textup{min}\left\{\mathcal{Q}^* \right\}\lambda_\textup{min}\left\{\mathcal{P}^*\right\}^\frac{p-1}{p}}{\lambda_\textup{max}\left\{\mathcal{P}^*\right\}} \frac{p}{3p-2}. 
			\end{split}
		\end{align}
		Therefore, based on the inequality \eqref{eqn:Lyapunov Derivative 3}, Lemma \ref{lem:FTS} and Lemma \ref{lem:FFTS}, we conclude that the origin of the error dynamics \eqref{eqn:Differentiator} is finite-time stable and fast finite-time stable.
	}
	\qedsymbol{}
\end{Proof*}
\begin{Corollary}[Disturbance Robustness]\label{cor:FTS Differentiator Disturbance Robustness}
	Consider the proposed HC-FFTSD \eqref{eqn:Differentiator} in Theorem \ref{thr:FFTS Differentiator} under perturbation, $ \delta = (\delta_1\Tp, \delta_2\Tp)\Tp$, $\delta_1, \delta_2 \in \bR^n$, and $\delta$ is bounded as $||\delta||\leq \bar{\delta}$. Thereafter, the differentiator under perturbation is as
	\begin{align}\label{eqn:Differentiator Perturbation}
		\begin{array}{ll}
			\dot{e}_1 &= -k_1 \phi_1(e_1) + e_2+\delta_1, \\ 
			\dot{e}_2 &= -k_2 \phi_2(e_1) + \delta_2.
		\end{array}    
	\end{align}
	When $\gamma_1$ in \eqref{eqn:gamma1 gamma2} fulfills $\gamma_1 \geq {\lambda_\textup{max}\left\{\mathcal{P}\right\}}/ {\lambda_\textup{min}\left\{\mathcal{P}\right\}}$, \eqref{eqn:Differentiator Perturbation} is practically finite-time stable (PFTS).
\end{Corollary}
\begin{Proof*}
	{\em
		Consider the Lyapunov stability analysis in Theorem \ref{thr:FFTS Differentiator}. With the Lyapunov-candidate defined by \eqref{eqn:Lyapunov FTS Differentiator} and the expression of the differentiator under perturbation in \eqref{eqn:Differentiator Perturbation}, we express the time derivative of \eqref{eqn:Lyapunov FTS Differentiator} as follows:
		\begin{align}\label{eqn:Lyapunov Derivative 4}
			\dot{V} & \leq  -\gamma_1 V -\gamma_2 V^\frac{1}{p} + 2\lambda_\textup{max}\left\{\mathcal{P}\right\}\bar{\delta} ||\zeta||.
		\end{align}
		By applying Cauchy-Schwarz inequality and \eqref{eqn:Lyapunov Bound 1}, from \eqref{eqn:Lyapunov Derivative 4}, we obtain,
		\begin{align}\label{eqn:Lyapunov Derivative 5}
			\begin{split}
				\dot{V} & \leq  -\gamma_1 V -\gamma_2 V^\frac{1}{p} + 
				\lambda_\textup{max}\left\{\mathcal{P}\right\} ||\zeta||^2 + \lambda_\textup{max}\left\{\mathcal{P}\right\} \bar{\delta}^2 \\
				& \leq -\left(\gamma_1- \frac{\lambda_\textup{max}\left\{\mathcal{P}\right\}}{\lambda_\textup{min}\left\{\mathcal{P}\right\}}\right) V - \gamma_2 V^\frac{1}{p} + \lambda_\textup{max}\left\{\mathcal{P}\right\} \bar{\delta}^2.
			\end{split}
		\end{align}
		Therefore, according to Lemma \ref{lem:PFTS}, with inequality \eqref{eqn:Lyapunov Derivative 5}, we conclude that the system \eqref{eqn:Differentiator Perturbation}, which is the differentiator \eqref{eqn:Differentiator} under disturbance $\delta$, is practical finite time stable at the origin.
	}\qedsymbol{}
\end{Proof*}
\begin{Corollary}[Noise Robustness]\label{cor:FTS Differentiator Noise Robustness}
	Consider the proposed HC-FFTSD \eqref{eqn:Differentiator} in Theorem \ref{thr:FFTS Differentiator} under measurement noise $\mu$, so that $\phi_1(e_1)$ and $\phi_2(e_1)$ in \eqref{eqn:phi1phi2} are replaced by $\phi_1(e_1+\mu)$ and $\phi_2(e_1+\mu)$ 
	in the differentiator, as follows: 
	\begin{align}\label{eqn:Differentiator Noise 1}
		\begin{array}{ll}
			\dot{e}_1 &= -k_1 \phi_1(e_1+\mu) + e_2 \\ 
			\dot{e}_2 &= -k_2 \phi_2(e_1+\mu), 
		\end{array}    
	\end{align}
	where $\mu$ is bounded as $||\mu||\leq \bar{\mu}$. When $\gamma_1$ in \eqref{eqn:gamma1 gamma2} fulfills $\gamma_1 \geq {\lambda_\textup{max}\left\{\mathcal{P}\right\}}/ {\lambda_\textup{min}\left\{\mathcal{P}\right\}}$, \eqref{eqn:Differentiator Noise 1} is practically finite-time stable (PFTS).
\end{Corollary}
\begin{Proof*}
	{\em
		From \eqref{eqn:Differentiator Noise 1}, we obtain the following expression
		\begin{align}\label{eqn:Differentiator Noise 2}
			\begin{array}{ll}
				&\dot{e}_1 = -k_1 \phi_1(e_1)+ e_2+k_1 \phi^*_1(e_1,\mu),  \\ 
				&\dot{e}_2 = -k_2 \phi_2(e_1)+k_2 \phi^*_2(e_1,\mu), \\
				&\phi^*_1(e_1,\mu) = -\phi_1(e_1+\mu) + \phi_1(e_1), \\
				&\phi^*_2(e_1,\mu) = -\phi_2(e_1+\mu) + \phi_2(e_1).
			\end{array}    
		\end{align}
		From \eqref{eqn:phi1phi2}, we obtain 
		\begin{align*}
			\phi^*_1(e_1,\mu) &= -\phi_1(e_1+\mu) + \phi_1(e_1) \\
			&= -k_3\mu - \|e_1+\mu\|^\frac{2(1-p)}{3p-2}(e_1+\mu) + \|e_1\|^\frac{2(1-p)}{3p-2}e_1 \\
			\phi^*_2(e_1,\mu) &= -\phi_2(e_1+\mu) + \phi_2(e_1)\\
			&= -k^2_3\mu - \frac{2k_3(2p-1)}{3p-2}\|e_1+\mu\|^\frac{2(1-p)}{3p-2}(e_1+\mu)-\frac{p}{3p-2}\|e_1+\mu\|^\frac{4(1-p)}{3p-2}(e_1+\mu)\\
			&+\frac{2k_3(2p-1)}{3p-2} \|e_1\|^\frac{2(1-p)}{3p-2}e_1+\frac{p}{3p-2}\|e_1\|^\frac{4(1-p)}{3p-2}e_1.
		\end{align*}
		Therefore, according to Lemma \ref{lem:Inequality Noise Robustness}, we obtain the upper bounds of $\|\phi^*_1(e_1,\mu)\| $ and $\|\phi^*_2(e_1,\mu)\| $ as:
		\begin{align*}
			\|\phi^*_1(e_1,\mu)\|& \leq k_3 \bar{\mu} + 2^{\frac{2(p-1)}{3p-2}} (\bar{\mu})^{1-\frac{2(p-1)}{3p-2}}\\
			\|\phi^*_2(e_1,\mu)\|& \leq k^2_3 \bar{\mu} +  \frac{2k_3(2p-1)}{3p-2}2^{\frac{2(p-1)}{3p-2}}(\bar{\mu})^{1-\frac{2(p-1)}{3p-2}} +\frac{p}{3p-2}2^{\frac{4(p-1)}{3p-2}}(\bar{\mu})^{1-\frac{4(p-1)}{3p-2}}.
		\end{align*}
		Applying Corollary \ref{cor:FTS Differentiator Disturbance Robustness} and with the above  bounds on $\|\phi^*_1(e_1,\mu)\|$ and $\|\phi^*_1(e_1,\mu)\|$, we conclude that the error dynamics \eqref{eqn:Differentiator Noise 1} is PFTS at the origin.
	}\qedsymbol{}
\end{Proof*}
\section{Problem Formulation}\label{sec:Problem}
\subsection{Coordinate frame definition}
\begin{figure}[ht]
	\centering
	\includegraphics[width=0.45\columnwidth]{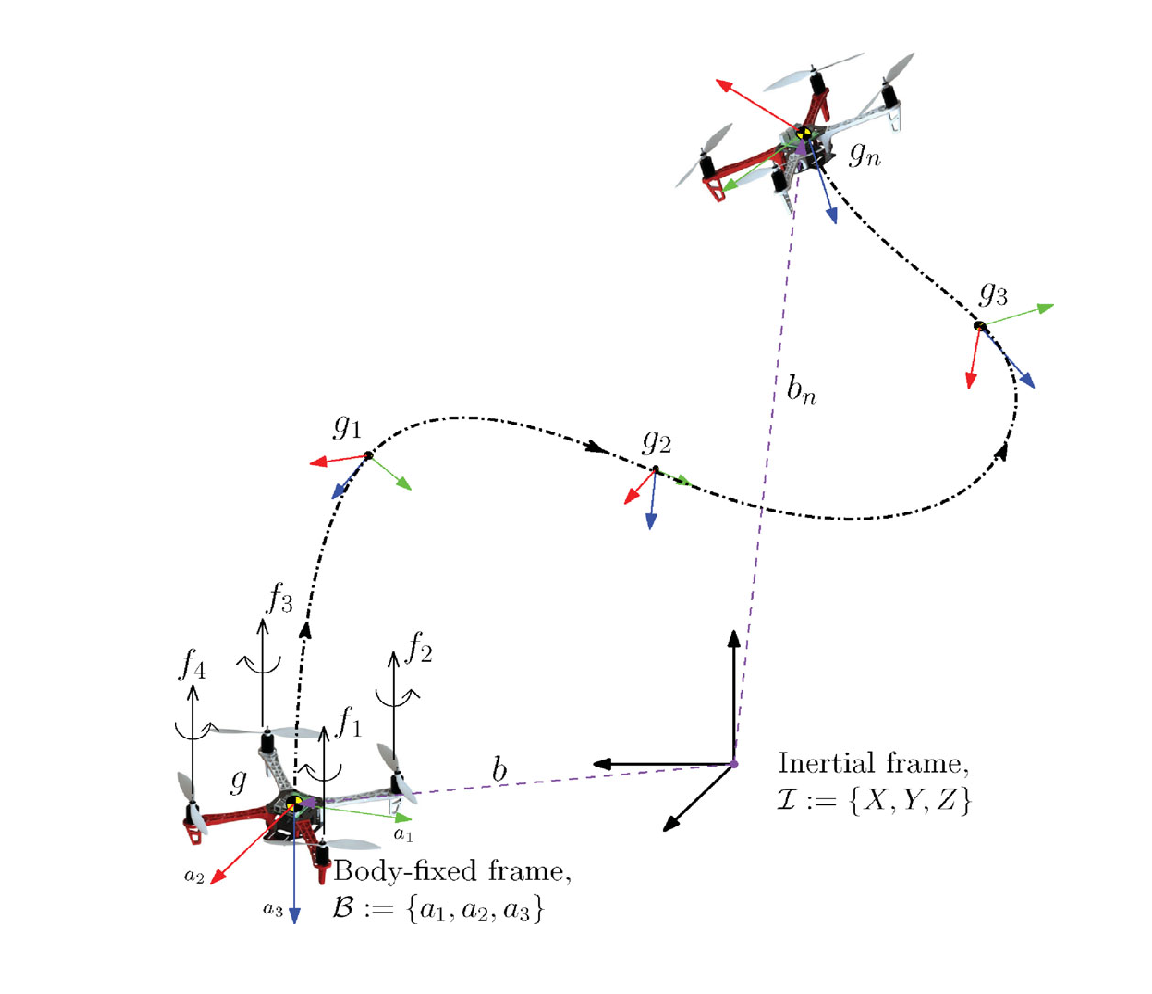}
	\caption{Guidance of a rotorcraft UAV through a trajectory between initial and final configurations on \SE \ \cite{hamrah2022finite,viswanathan2018integrated}}
	\label{fig:Configuration}
\end{figure}

The configuration of the vehicle, modeled as a rigid body, is given by its position and orientation, which are together
referred to as its pose. To define the pose of the vehicle, we fix a coordinate frame $\mathcal{B}$ to its body and another coordinate frame $\mathcal{I}$ that is fixed in space as the inertial coordinate frame. Define $\textbf{e}_i$ as the unit vector along the $i$th coordinate axis for $i=1,2,3$. Let $b \in \bR^3$ denote the position vector of the origin of frame $\mathcal{B}$ with respect to frame $\mathcal{I}$. Let \SO \ denote the orientation (attitude), defined as the rotation matrix from frame  $\mathcal{B}$ to frame  $\mathcal{I}$. The pose of the vehicle can be represented in matrix form as follows:
\begin{align}\label{eqn:Coordinate 1}
	g = 
	\left[
	\begin{array}{cc}
		R & b \\
		0 & 1 \\
	\end{array}
	\right]
	\in \SE
\end{align}
where \SE, the special Euclidean group, is the six-dimensional Lie group of rigid body motions. A diagram of guidance and trajectory tracking on \SE \ through a set of waypoints is presented in Figure \ref{fig:Configuration} as follows.
\subsection{System kinematics and dynamics}
The instantaneous pose (position and attitude) is compactly represented by $\mathrm{g} = (b, R) \in  \SE$. The UAV's kinematics is then defined by:
\begin{align}
	\left\{
	\begin{array}{lr}
		\dot{b} = v = R\nu,\\
		\dot{R} =R\Omega^{\times},
	\end{array}
	\right.
	\label{eqn:Kinematics}
\end{align}
where $v \in \R{3}$ and $\nu \in \R{3}$ denote the translational velocity in frames $\mathcal{I}$ and $\mathcal{B}$ respectively, and $\Omega$ $\in$ $ \R{3}$ is the angular velocity in body-fixed frame $\mathcal{B}$. The overall system kinematics and dynamics of a rotorcraft UAV with a body-fixed plane of rotors are given by: 
\begin{align}\label{eqn:System}
	\begin{cases}
		& \dot{b} = v = R\nu \\
		& m\dot{v} = mg\textbf{e}_3 - fR\textbf{e}_3 + \varphi_D  \\
		& \dot{R} = R\Omega^{\times} \\
		&J\dot{\Omega} = J\Omega\times\Omega + \tau + \tau_D 
	\end{cases}    
\end{align}
where $\textbf{e}_3 = [0 \quad 0 \quad 1]\Tp$, $f \in \bR$ is the scalar thrust force, and $\tau \in \bR^{3}$ is the control torque created by the rotors, $g$ denotes the acceleration due to gravity and $m \in \bR^{+}$ and $J = J\T \in \bR^{3 \times 3}$ are the mass and inertia matrix of the UAV, respectively. The force and torque disturbances are denoted $\varphi_D$ and $\tau_D$ respectively, which are mainly due to unsteady aerodynamics.	
\subsection{Morse function on $\SO$}
The following Lemma is used for the rotational ESO scheme designed for the vehicle.
\begin{Lemma}\label{lem:sK Definition}
	\textup{\cite{bohn2016almost}}
	Consider attitude kinematics 
	\begin{align}\label{eqn:Attitude Kinematics}
		\dot{R}=R\Omega^\times, R\in\SO, \Omega\in \sso.
	\end{align}
	Define $K= \diag([K_1,K_2,K_3])$, where $K_1>K_2>K_3 \ge 1$. Define
	\begin{align}\label{eqn:sK Definition}
		s_K(R)=\sum\limits_{i=1}^{3}K_i(R\Tp \textbf{\textup e}_i)\times \textbf{\textup e}_i, 
	\end{align}
	such that\textup{ $\frac{\di}{\di t}\lan K,I-R\ran=\Omega\T s_K(R)$}. Here \textup{$\lan A, B\ran=\tr  (A\T B)$}, which makes 
	$\lan K,I-R\ran$ a Morse function defined on $\SO$.
	Let $\cS\subset\SO$ be a closed subset containing the identity in its interior, defined by
	\begin{align}\label{eqn:S Definition}
		\begin{split}
			\cS &= \big\{ R\in\SO\, :\, R_{ii}\ge 0 \mbox{\textup{ and} } R_{ij}R_{ji}\le 0, \forall i,~ j\in \{1,2,3\},\ i\ne j\big\}. 
		\end{split}
	\end{align}
	Then for $\forall R\in\cS$, we have 
	\begin{align}\label{eqn:sK Bound} 
		s_K(R)\Tp s_K(R) \geq \langle K, I-R \rangle.  
	\end{align}
\end{Lemma}
\begin{Remark}[Almost global domain of attraction]\label{rem:Critical Points}
	{\em 
		\cite{sanyal2010almost} We know that the subset of $\SO$ where $s_K(R) = 0, R \in \SO $, which is also the set of critical points for $\langle I-R, K\rangle$, is 
		\begin{align}\label{eqn:Critical Points}
			\begin{split}
				C &\triangleq \{I, \textup{diag}(1,-1,-1), \textup{diag}(-1,1,-1), \textup{diag}(-1,-1,1)\}\subset \SO.
			\end{split}
		\end{align}
		In addition, the global minimum of this Morse function is $R=I$.
	}
\end{Remark}
\subsection{ESO estimates and errors}
The ESO design on $\SE$ is split into a translational ESO design on vector space $\bR^3$ and a rotational ESO design on $\SO$. Let ($\widehat{b}, \widehat{v}, \widehat{\varphi}_D )\in \bR^3 \times \bR^3 \times \bR^3$ be the estimated position, translational velocity, and disturbance force, as the states of the translational ESO. The estimation errors for the translational ESO are defined as follows:
\begin{align}\label{eqn:Translational ESO Error}
	e_b = b - \widehat{b}, e_v = v -\widehat{v}, e_\varphi = \varphi_D - \widehat{\varphi}_D,
\end{align}
which are estimation errors of position, translational velocity, and total disturbance force respectively. 

Let $(\widehat{R}, \widehat{\Omega}, \widehat{\tau}_D)\in \SO\times\bR^3\times\bR^3$ be the estimated attitude, angular velocity, and disturbance torque states provided by the rotational ESO. For the rotational ESO, the error states are defined as follows. The attitude estimation error is defined as:
\begin{align}\label{eqn:Attitude ESO Error ER}
	E_R = \widehat{R}\Tp R,
\end{align}
on the group of rigid body rotations, $\SO$,  which is not a vector space. 
The angular velocity estimation error, $e_\Omega$, and torque disturbance estimation error, $e_\tau$, are expressed on the vector space $\bR^3$, and are defined as:
\begin{align}\label{eqn:Attitude ESO Error}
	e_\Omega = \Omega -E\Tp_R\widehat{\Omega}, \ e_\tau = \tau_D - \widehat{\tau}_D.
\end{align}

With a proper ESO design on $\SE$, the error states $(e_b, e_v, e\varphi)$ and $(E_R, e_\Omega, e_\tau)$ should converge to $(0, 0, 0)$ and $(I, 0, 0)$, respectively. The ESO design and its stability proof are described in detail in the following section.

\section{Fast Finite-Time Stable Extended State Observer (FFTS-ESO) on \SE}\label{sec:ESO}
In this section, we present the FFTS-ESO on $\SE$. As mentioned in the previous section, the ESO design on $\SE$ can be represented as a translational ESO on the vector space $\bR^3$ to estimate disturbance forces, and an rotational ESO on $\SO$ to estimate disturbance torques. We present the two ESO designs along their stability proofs in this section.  
\subsection{ESO for Translational Motion}
\begin{Proposition}[Translational ESO]\label{prop:Translational ESO}
	Define the positive scalar gains $k_{t1}$ and $k_{t2}$, which make the matrix $\mathcal{A}_t \in \bR^{2\times2}$ defined as:
	\begin{align}\label{eqn:Matrix At}
		\mathcal{A}_t = 
		\begin{bmatrix}
			-k_{t1} & 1 \\ 
			-k_{t2} & 0
		\end{bmatrix},
	\end{align}
	a Hurwitz matrix. The ESO designed for the translational motion is given by:
	\begin{align}\label{eqn:Translational ESO}
		\begin{split}
			\dot{\widehat{b}} &= \widehat{v}, \\
			m\dot{\widehat{v}} &= mg\textbf{\textup{e}}_3 -fR\textbf{\textup{e}}_3 + m k_{t1} \phi_1(\psi_t) + m\kappa_t \Big[ (e\Tp_be_b)^{\frac{1-p}{p}}H\Big(e_b,\frac{p-1}{p}\Big)e_v + e_v\Big] +\widehat{\varphi}_D , \\
			\dot{\widehat{\varphi}}_D &= m k_{t2} \phi_2(\psi_t),
		\end{split}
	\end{align}
	where $\psi_t$ is defined as
	\begin{align}\label{eqn:psit}
		\psi_t = e_v + \kappa_t \Big[e_b + (e\Tp_be_b)^{\frac{1-p}{p}}e_b\Big],\; \kappa_t>1/2,
	\end{align}
	and $\phi_1(\cdot)$ is as defined by the expression in \eqref{eqn:phi1phi2}. In addition, the constant $k_{t3}$ is defined 
	and it occurs in the terms $\phi_1(\psi_t)$ and $\phi_2(\psi_t)$, where it takes the place of $k_3$ in \eqref{eqn:phi1phi2}.
\end{Proposition}
\begin{Theorem}\label{thr:Transaltional ESO Error Dynamics}
	With the observer errors for the translational ESO defined by \eqref{eqn:Translational ESO Error}, the translational kinematics and dynamics given by \eqref{eqn:System}, and the ESO for translational motion given in Proposition \ref{prop:Translational ESO}, the error dynamics of the ESO is given by:
	\begin{align}\label{eqn:Translational ESO Error Dynamics}
		\begin{split}
			\dot{e}_b &= e_v, \\
			m\dot{e}_v &= -mk_{t1} \phi_1(\psi_t)  - m\kappa_t \left[ (e\Tp_be_b)^{\frac{1-p}{p}}H\left(e_b,\frac{p-1}{p}\right)e_v + e_v\right] + e_\varphi,\\
			\dot{e}_\varphi &= -mk_{t2}\phi_2(\psi_t)  + \dot{\varphi}_D.
		\end{split}
	\end{align}
	The error dynamics 
	\eqref{eqn:Translational ESO Error Dynamics} is FFTS at the 
	origin $((e_b,e_v,e_\varphi)=(0,0,0))$, when the resultant disturbance force is constant $(\dot{\varphi}_D=0)$ and the observer gains are constrained according to Proposition \ref{prop:Translational ESO}.
\end{Theorem}
\begin{Proof*}
	{\em
		Simplify \eqref{eqn:Translational ESO Error Dynamics} as:
		\begin{align}\label{eqn:Translational ESO Error Dynamics 1}
			\begin{split}
				\dot{\psi}_t &= -k_{t1} \phi_1(\psi_t) + m^{-1} e_\varphi,  \\
				m^{-1}\dot{e}_\varphi &= -k_{t2}\phi_2(\psi_t)  + m^{-1} \dot{\varphi}_D. 
			\end{split}
		\end{align}
		Next, for $\mathcal{A}_t$ as defined in \eqref{eqn:Matrix At}, $\forall\, \mathcal{Q}_t \in \bR^{ 2\times 2}$ where $\mathcal{Q}_t \succ 0$, the Lyapunov equation,
		\begin{align}\label{eqn:Lyapunov Equation TESO 1}
			\mathcal{A}\Tp_t \mathcal{P}_t+\mathcal{P}_t\mathcal{A}_t =-\mathcal{Q}_t,
		\end{align}
		has a unique solution $\mathcal{P}_t$. Thereafter, define the Lyapunov function: 
		\begin{align}\label{eqn:Lyapunov Translational ESO}
			V_t &= V_{t0}  + \mu_te\Tp_be_b, \mbox{ where } V_{t0} = \zeta\Tp_t \mathcal{P}_t \zeta_t 
		\end{align}
		and $\zeta_t$ is defined as:
		\begin{align*}
			\zeta_t = [\phi\Tp_1(\psi_t), \ m^{-1} e\Tp_\varphi]\Tp.
		\end{align*}
		We constrain the positive scalar $\mu_t$ in \eqref{eqn:Lyapunov Translational ESO} as:
		\begin{align}\label{eqn:Translational ESO Constraint}
			0<\mu_t<k^3_{t3}\frac{\lambda_\textup{min}\left\{\mathcal{P}_t\right\}\lambda_\textup{min}\left\{\mathcal{Q}_t\right\}}{ \lambda_\textup{max}\left\{\mathcal{P}_t\right\}}.
		\end{align}
		From Theorem \ref{thr:FFTS Differentiator}, \eqref{eqn:Translational ESO Error Dynamics 1} and \eqref{eqn:Lyapunov Bound 3}, we find that the time-derivative of $V_t$ satisfies:
		\begin{align}\label{eqn:Lyapunov Derivative Translational ESO 1}
			\dot{V}_t &\leq -\gamma_{t1} V_{t0} -\gamma_{t2} V^\frac{1}{p}_{t0} + 2\mu_te\Tp_be_v, 
		\end{align}
		where $\gamma_{t1}$ and $\gamma_{t2}$ are defined by:
		\begin{align}\label{eqn:gammat1 gammat2}
			\begin{split}
				\gamma_{t1} = k_{t3} \frac{\lambda_\textup{min}\left\{\mathcal{Q}_t \right\}}{\lambda_\textup{max}\left\{\mathcal{P}_t\right\}},\; 
				\gamma_{t2} = \frac{\lambda_\textup{min}\left\{\mathcal{Q}_t \right\}\lambda_\textup{min}\left\{\mathcal{P}_t\right\}^\frac{p-1}{p}p}{\lambda_\textup{max}\left\{\mathcal{P}_t\right\}(3p-2)}.
			\end{split}
		\end{align}
		Substituting \eqref{eqn:psit} into \eqref{eqn:Lyapunov Derivative Translational ESO 1}, we obtain:
		\begin{align}\label{eqn:Lyapunov Derivative Translational ESO 2}
			\begin{split}
				\dot{V}_t &\leq -\gamma_{t1} V_{t0} -\gamma_{t2} V^\frac{1}{p}_{t0} +2 \mu_t e\Tp_b\Big[\psi_t - \kappa_t e_b- \kappa_t (e\Tp_be_b)^{\frac{1-p}{p}}e_b\Big] \\
				&\leq -\gamma_{t1} V_{t0} -\gamma_{t2} V^\frac{1}{p}_{t0} +2\mu_t  e\Tp_b \psi_t - 2 \mu_t\kappa_t e\Tp_be_b- 2 \mu_t\kappa_t (e\Tp_be_b)^{\frac{1}{p}}  \\
				&\leq -\gamma_{t1} V_{t0} -\gamma_{t2} V^\frac{1}{p}_{t0} - 2 \mu_t\kappa_t e\Tp_be_b- 2\mu_t \kappa_t (e\Tp_be_b)^{\frac{1}{p}} + \mu_t \psi\Tp_t \psi_t + \mu_t e\Tp_be_b \\
				&\leq -\left(\gamma_{t1}-\frac{\mu_t}{k^2_{t3}\lambda_\textup{min}\left\{\mathcal{P}_t\right\}}\right) V_{t0} -\gamma_{t2} V^\frac{1}{p}_{t0} - (2\kappa_t-1) \mu_t e\Tp_be_b- 2 \kappa_t \mu_t^\frac{p-1}{p}\mu_t^\frac{1}{p} (e\Tp_be_b)^{\frac{1}{p}}. 
			\end{split}
		\end{align}
		Therefore, we further obtain:
		\begin{align}\label{eqn:Lyapunov Derivative Translational ESO 3}
			\dot{V}_t < -\Gamma_{t1} V_t  - \Gamma_{t2} V^{\frac{1}{p}}_t,
		\end{align}
		where 
		\begin{align}\label{eqn:Gammat1 Gammat2}
			\begin{split}
				\Gamma_{t1}  &= \textup{min}\left\{k_{t3} \frac{\lambda_\textup{min}\left\{\mathcal{Q}_t \right\}}{\lambda_\textup{max}\left\{\mathcal{P}_t\right\}} - \frac{\mu_t}{k^2_{t3} \lambda_\textup{min}\left\{\mathcal{P}_t\right\}}, 2\kappa_t - 1 \right\}, \\
				\Gamma_{t2}  &=  \textup{min}\left\{\frac{\lambda_\textup{min}\left\{\mathcal{Q}_t \right\}\lambda_\textup{min}\left\{\mathcal{P}_t\right\}^\frac{p-1}{p}p}{\lambda_\textup{max}\left\{\mathcal{P}_t\right\}(3p-2)}, 2\kappa_t\mu_t^\frac{p-1}{p} \right\} .
			\end{split}
		\end{align}
		Based on \eqref{eqn:Lyapunov Derivative Translational ESO 3}, we conclude that when the resultant disturbance force is constant and the ESO gains satisfy the constraints 1-3 in Proposition \ref{prop:Translational ESO}, the error dynamics of the ESO  \eqref{eqn:Translational ESO Error Dynamics} is FFTS. This concludes the proof of Theorem \ref{thr:Transaltional ESO Error Dynamics}.
	}
	\qedsymbol{}
\end{Proof*}
\subsection{ESO for Rotational Motion}
\begin{Proposition}[Rotational ESO]\label{prop:Attitude ESO}
	Define $e_R = s_K(E_R)$, where $s_K(\cdot)$ is as defined by Lemma \ref{lem:sK Definition}. 
	Define $e_w(E_R,e_\Omega)$ as follows:
	\begin{align}\label{eqn:ew}
		e_w(E_R,e_\Omega) = \frac{\textup{d}}{\textup{d}t} e_R =  \sum_{i=1}^3K_i\textbf{\textup{e}}_i \times (e_\Omega\times E\Tp_R \textbf{\textup{e}}_i).
	\end{align}
	Define the positive scalar gains $k_{a1}$ and $k_{a2}$, which make the matrix $\mathcal{A}_a \in \bR^{2\times2}$ defined as:
	\begin{align}\label{eqn:Matrix Aa}
		\mathcal{A}_a = 
		\begin{bmatrix}
			-k_{a1} & 1 \\ 
			-k_{a2} & 0
		\end{bmatrix},
	\end{align}
	a Hurwitz matrix. The ESO designed for the rotational motion is given by:
	\begin{align}\label{eqn:Attitude ESO}
		\begin{split}
			&\dot{\widehat{R}}= \widehat{R}\widehat{\Omega}^\times, \\
			&\dot{\widehat{\Omega}}= E_RJ^{-1}\left[J\Omega\times\Omega+\widehat{\tau}_D +\tau + k_{a1} J\phi_1(\psi_a)+\kappa_a J (e\Tp_Re_R)^{\frac{1-p}{p}}H\Big(e_R,\frac{p-1}{p}\Big)e_w \right]\\
			& + E_RJ^{-1}(\kappa_a Je_w)+E_Re^\times_\Omega E\Tp_R\widehat{\Omega}, \\
			&\dot{\widehat{\tau}}_D = J k_{a2} \phi_2(\psi_a), 
		\end{split}
	\end{align}
	where $\psi_a$ is defined as follows: 
	\begin{align}\label{eqn:psia}
		\psi_a = e_\Omega + \kappa_a \Big[e_R + (e_R\Tp e_R)^{\frac{1-p}{p}}e_R\Big], \ \kappa_a>\frac{1}{2}.
	\end{align}  
	In addition, the constant $k_{a3}$ is defined and it occurs in the terms $\phi_1(\psi_a)$ and $\phi_2(\psi_a)$,
	where it takes the place of $k_3$ in \eqref{eqn:phi1phi2}. 
\end{Proposition}
\begin{Theorem}\label{thr:Attitude ESO Error Dynamics}
	With the observer errors for the rotational ESO defined by \eqref{eqn:Attitude ESO Error}, the rotational kinematics and dynamics given by \eqref{eqn:System}, and the ESO for rotational motion given in Proposition \ref{prop:Attitude ESO}, the error dynamics of the ESO is given by:
	\begin{align}\label{eqn:Attitude ESO Error Dynamics}
		\begin{split}
			\dot{E}_R&= E_Re^\times_\Omega, \\
			J\dot{e}_\Omega &= -k_{a1}J \phi_1(\psi_a) - \kappa_a J\left[(e\Tp_Re_R)^{\frac{1-p}{p}}H\Big(e_R,\frac{p-1}{p}\Big)e_w + e_w\right] + e_\tau,  \\
			\dot{e}_\tau &= -k_{a2}J\phi_2(\psi_a) +\dot{\tau}_D. 
		\end{split}
	\end{align}
	The error dynamics \eqref{eqn:Attitude ESO Error Dynamics} is almost globally FFTS at the origin $((E_R,e_\Omega,e_\tau)=(I,0,0))$, when the resultant disturbance torque is constant $(\dot{\tau}_D=0)$ and the observer gains are constrained  according to Proposition \ref{prop:Attitude ESO}. 
\end{Theorem}
\begin{Proof*}
	{\em
		Simplify \eqref{eqn:Attitude ESO Error Dynamics} as:
		\begin{align}\label{eqn:Attitude ESO Error Dynamics 1}
			\begin{split}
				\dot{\psi}_a &= -k_{a1} \phi_1(\psi_a) + J^{-1} e_\tau, \\
				J^{-1}\dot{e}_\tau &= -k_{a2}\phi_2(\psi_a)+J^{-1}\dot{\tau}_D. 
			\end{split} 
		\end{align}
		Next, for $\mathcal{A}_a$ as defined in \eqref{eqn:Matrix Aa}, 
		$\forall\, \mathcal{Q}_a \in \bR^{ 2\times 2}$ where $\mathcal{Q}_a  \succ 0$, the Lyapunov equation: 
		\begin{align}\label{eqn:Lyapunov Equation RESO 1}
			\mathcal{A}\Tp_a \mathcal{P}_a+\mathcal{P}_a\mathcal{A}_a =-\mathcal{Q}_a,
		\end{align}
		has a unique solution $\mathcal{P}_a$.	
		Thereafter, define the Morse-Lyapunov function: 
		\begin{align}\label{eqn:Lyapunov Attitude ESO}
			\begin{split}
				V_a &= V_{a0}  + \mu_a\langle K, I-E_R \rangle, \mbox{ where } V_{a0} = \zeta\Tp_a \mathcal{P}_a \zeta_a,
			\end{split}
		\end{align}
		$\mu_a$ is a positive scalar, and $\zeta_a$ is defined as:
		\begin{align*}
			\zeta_a = [\phi\Tp_1(\psi_a), \ J^{-1} e\Tp_\tau]\Tp.
		\end{align*}
		We constrain the positive scalar $\mu_a$ in \eqref{eqn:Lyapunov Attitude ESO} as:
		\begin{align}\label{eqn:Attitude ESO Constraint}
			0<\mu_a<2k^3_{a3}\frac{\lambda_\textup{min}\left\{\mathcal{P}_a\right\}\lambda_\textup{min}\left\{\mathcal{Q}_a\right\}}{ \lambda_\textup{max}\left\{\mathcal{P}_a\right\}}.
		\end{align}
		From Theorem \ref{thr:FFTS Differentiator}, \eqref{eqn:Attitude ESO Error Dynamics 1} and \eqref{eqn:Lyapunov Bound 3}, we find that the time-derivative of $V_a$ satisfies:
		\begin{align}\label{eqn:Lyapunov Derivative Attitude ESO 1}
			\dot{V}_a &\leq -\gamma_{a1} V_{a0} -\gamma_{a2} V^\frac{1}{p}_{a0} + \mu_a e\Tp_Re_\Omega,
		\end{align}
		where $\gamma_{a1}$ and $\gamma_{a2}$ are defined by:
		\begin{align}\label{eqn:gammaa1 gammaa2}
			\begin{split}
				\gamma_{a1} = k_{a3} \frac{\lambda_\textup{min}\left\{\mathcal{Q}_a \right\}}{\lambda_\textup{max}\left\{\mathcal{P}_a\right\}}, \;
				\gamma_{a2} = \frac{\lambda_\textup{min}\left\{\mathcal{Q}_a \right\}\lambda_\textup{min}\left\{\mathcal{P}_a\right\}^\frac{p-1}{p}p}{\lambda_\textup{max}\left\{\mathcal{P}_a\right\}(3p-2)}.
			\end{split}
		\end{align}
		Substituting \eqref{eqn:psia} into \eqref{eqn:Lyapunov Derivative Attitude ESO 1}, we obtain:
		\begin{align}\label{eqn:Lyapunov Derivative Attitude ESO 2}
			\begin{split}
				\dot{V}_a &\leq -\gamma_{a1} V_{a0} -\gamma_{a2} V^\frac{1}{p}_{a0} + \mu_a e\Tp_R\Big[\psi_a - \kappa_a e_R- \kappa_a (e\Tp_Re_R)^{\frac{1-p}{p}}e_R\Big] \\
				&\leq -\gamma_{a1} V_{a0} -\gamma_{a2} V^\frac{1}{p}_{a0} + \frac{1}{2}\mu_a \left(e\Tp_R e_R +\psi\Tp_a\psi_a \right)- \kappa_a \mu_a \left[e\Tp_Re_R+ (e\Tp_Re_R)^{\frac{1}{p}} \right] \\
				&\leq -\left(\gamma_{a1}-\frac{\mu_a}{2k^2_{a3}\lambda_\textup{min}\left\{\mathcal{P}_a\right\}}\right) V_{a0} -\gamma_{a2} V^\frac{1}{p}_{a0} -\left(\kappa_a-\frac{1}{2}\right) \mu_a e\T_Re_R- \kappa_a \mu_a (e\T_Re_R)^{\frac{1}{p}}.
			\end{split}
		\end{align}
		By applying Lemma \ref{lem:sK Definition} on \eqref{eqn:Lyapunov Derivative Attitude ESO 1}, we obtain:	
		\begin{align}\label{eqn:Lyapunov Derivative Attitude ESO 3}
			\begin{split}
				\dot{V}_a &\leq -\left(\gamma_{a1}-\frac{\mu_a}{2k^2_{a3}\lambda_\textup{min}\left\{\mathcal{P}_a\right\}}\right) V_{a0} -\gamma_{a2} V^\frac{1}{p}_{a0} \\
				&- \left(\kappa_a-\frac{1}{2}\right) \mu_a \langle K,I-E_R\rangle- \kappa_a \mu^\frac{p-1}{p}_a \mu^\frac{1}{p}_a \langle K,I-E_R\rangle^{\frac{1}{p}}.
			\end{split}
		\end{align}	
		After some algebra, we further obtain:
		\begin{align}\label{eqn:Lyapunov Derivative Attitude ESO 4}
			\dot{V}_a \leq -\Gamma_{a1} V_a  - \Gamma_{a2} V^{\frac{1}{p}}_a,
		\end{align}		
		where: 
		\begin{align}\label{eqn:Gammaa1 Gammaa2}
			\begin{split}
				\Gamma_{a1}&=\textup{min}\left\{k_{a3}\frac{\lambda_\textup{min}\left\{\mathcal{Q}_a \right\}}{\lambda_\textup{max}\left\{\mathcal{P}_a\right\}}-\frac{\mu_a}{2k^2_{a3}\lambda_\textup{min}\left\{\mathcal{P}_a\right\}},\kappa_a-\frac{1}{2}\right\},\\
				\Gamma_{a2}&=\textup{min}\left\{\frac{\lambda_\textup{min}\left\{\mathcal{Q}_a \right\}\lambda_\textup{min}\left\{\mathcal{P}_a\right\}^\frac{p-1}{p}p}{\lambda_\textup{max}\left\{\mathcal{P}_a\right\}(3p-2)}, \kappa_a \mu^\frac{p-1}{p}_a \right\}.
			\end{split}
		\end{align}
		Considering the expression given by \eqref{eqn:Lyapunov Derivative Attitude ESO 4}, the set where $\dot{V}_a=0$ is:
		\begin{align}\label{eqn:Attitude ESO Equilibrium Point Set 1}
			\begin{split}
				&\dot{V}_a^{-1}(0)=\left\{(E_R,e_\Omega,e_\tau): s_K(E_R)=0, \mbox{and} \ \zeta_a = 0 \right\} \\
				&= \left\{(E_R,e_\Omega,e_\tau): E_R \in C, e_\Omega=0, \mbox{and} \ e_\tau=0 \right\},
			\end{split}
		\end{align}
		where $C$ is as defined by \eqref{eqn:Critical Points}, which gives the set of the critical points of the Morse function used as part of the Morse-Lyapunov function in \eqref{eqn:Lyapunov Attitude ESO}. Using Theorem 8.4 from \cite{khalil2002nonlinear}, we conclude that $(E_R,e_\Omega,e_\tau)$ converge to the set:
		\begin{align}\label{eqn:Attitude ESO Equilibrium Point Set 2}
			\begin{split}
				&S=\left\{(E_R,e_\Omega,e_\tau)\in \SO\times\bR^3\times\bR^3:  E_R \in C, e_\Omega=0, \mbox{and} \ e_\tau=0 \right\},
			\end{split}
		\end{align}
		in finite time.	
		Based on \eqref{eqn:Lyapunov Derivative Attitude ESO 4}, and Lemma \ref{lem:FFTS}, we conclude that when the observer gains satisfy the constraints in Proposition \ref{prop:Attitude ESO},  
		the error dynamics \eqref{eqn:Attitude ESO Error Dynamics} converges to the set $S$ in finite time. 
		
		In $S$, the only stable equilibrium is $(I,0,0)$, while the other three are unstable. The resulting closed-loop system with the estimation errors gives rise to a H\"{o}lder-continuous feedback with exponent less than one $(1/2<1/p<1)$, while in the limiting case of $p=1$, the feedback system is Lipschitz-continuous. Proceeding with a topological equivalence-based analysis similar to the one by \cite{bohn2016almost}, we conclude that the equilibrium and the corresponding regions of attraction of the rotational ESO with $p\in]1,2[$ are identical to those of the corresponding Lipschitz-continuous asymptotically stable ESO with $p=1$, and the region of attraction is almost global. 
		
		To summarize, we conclude that the error dynamics \eqref{eqn:Translational ESO Error Dynamics} is almost globally FFTS (AG-FFTS) at the origin $((E_R,e_\Omega,e_\tau)=(I,0,0))$ when the resultant disturbance torque is constant ($\dot{\tau}_D=0$) and the observer gains are constrained  according to Proposition \ref{prop:Attitude ESO}. This concludes the proof of Theorem \ref{thr:Attitude ESO Error Dynamics}. 
	}
	\qedsymbol{}
\end{Proof*}
\begin{Remark}[Disturbance robustness of the ESO]
	{\em
		Consider Corollary \ref{cor:FTS Differentiator Disturbance Robustness} and its constraints on differentiator gains. When the disturbance forces and torques are time-varying, then $\|\dot{\varphi}_D\|, \|\dot{\tau}_D\|>0$. Further, if the constraints on gains in Corollary \ref{cor:FTS Differentiator Disturbance Robustness} are fulfilled, the estimation error dynamics of the proposed ESO will be PFTS. 
	}
\end{Remark}
\begin{Remark}[Noise robustness of the ESO]\label{rem:Noise}
	{\em
		Consider Corollary \ref{cor:FTS Differentiator Noise Robustness} and its constraints on differentiator gains. When the ESO measurements have noise and the constraints on gains in Corollary \ref{cor:FTS Differentiator Noise Robustness} are fulfilled, the estimation error dynamics of the proposed ESO will be PFTS. Moreover, according to Lemma \ref{lem:PFTS} and Corollary \ref{cor:FTS Differentiator Noise Robustness}, the $\eta$ in \eqref{eqn:Preliminary PFTS Lyapunov Inequality} of Lemma \ref{lem:PFTS} is a function on the level of noise in information on $R$, $\Omega$, $b$ and $v$ and is monotonically increasing with the level of noise.
	}
\end{Remark}
\begin{Remark}[Comparative Analysis of Noise Robustness: FFTS-ESO vs. the FxTSDO by Liu et al.~\cite{liu2022fixed} ]\label{rem:Compare}
	{\em
		We investigate the disturbance (forces or torques) observers proposed by \cite{liu2022fixed} in their Theorems 1 and 2, known as FxTSDO. The input of FxTSDO relies on the motion signals, $X_2$, $Y_2$,  which represent translational and angular velocities, and $\dot{X}_2$, $\dot{Y}_2$, which represent translational and angular accelerations, respectively. 	
		However, the high-level noise associated with the translational acceleration obtained from an accelerometer restricts its direct use in a flight control scheme. Additionally, direct measurement of angular acceleration is usually not feasible. 
		Furthermore, if $\dot{X}_2$ and $\dot{Y}_2$ are obtained from the finite difference of $X_2$ and $Y_2$, they will have higher noise levels than $X_2$ and $Y_2$, leading to inferior disturbance estimation performance. 	
		In contrast to FxTSDO, the proposed FFTS-ESO incorporates position and attitude signals, which are zero-order derivatives of motions with lower noise levels. Consequently, FFTS-ESO outperforms FxTSDO in terms of disturbance estimation performance, despite the theoretical fixed-time stability of FxTSDO. We show this through our numerical simulations in Section \ref{sec:Numerical}.
	}
\end{Remark}

\section{Numerical Simulations}\label{sec:Numerical}
In this section, we compare the proposed FFTS-ESO with existing disturbance estimation schemes, which are LESO by \cite{shao2018robust} and FxTSDO by \cite{liu2022fixed}, on their disturbance estimation performance in four different simulated flight scenarios, with and without the presence of measurement noises. The four flight scenarios correspond to four desired trajectories. The inertia and mass of the simulated rotorcraft UAV are $  J=\textup{diag}([0.0820,0.0845,0.1377]) \ \textup{kg}\cdot \textup{m}^2, \quad m = 4.34 \ \textup{kg}$ by \cite{pounds2010modelling}. Since the goal of the simulation is to validate and compare the disturbance estimation performance, the actuator dynamics and saturation are not included in the results reported in this section. The tracking control scheme to drive the UAV to track the desired trajectories is developed based on the control scheme reported by \cite{viswanathan2017finite,viswanathan2018integrated}. As this scheme is not a contribution of this article, we omit its description for brevity. We use MATLAB/Simulink with its ODE2 (Heun method) solver to conduct this set of simulations. The time step size is $h = 0.001 \text{s}$ and the simulated duration is $T = 30 \text{s}$. 
\begin{table}[htbp]
	\begin{center}
		\begin{tabular}{ | c | c| } 
			\hline
			Hovering& $b_d(t)=\left[ 0, \ 0, \ -3\right]\Tp \textup{(m)}$  \\ 
			Slow Swing & $b_d(t)=\left[ 10 \ \textup{sin}(0.1\pi t), \ 0, \ -3\right]\Tp \textup{(m)}$ \\ 
			Fast Swing& $b_d(t)=\left[ 5 \ \textup{sin}(0.5\pi t), \ 0, \ -3\right]\Tp \textup{(m)}$ \\ 
			High Pitch& $b_d(t)=\left[ 10 \ \textup{sin}(0.5\pi t), \ 10 \ \textup{cos}(0.5\pi t), \ -3\right]\Tp \textup{(m)}$ \\
			\hline
		\end{tabular}
	\end{center}
	\caption{Flight trajectories to be tracked for the comparisons between  LESO, FxTSDO and FFTS-ESO}
	\label{table:ESO Comparison Trajectory}
\end{table}
\begin{table}[htbp]
	\begin{center}
		\begin{tabular}{ | c| c| c| } 
			\hline
			$b_\textup{N}$ & $b_\textup{N} = b + \mu_b$ & $ \mu_b \sim P_b = 3e^{-8} $   \\ 
			$v_\textup{N}$ & $v_\textup{N} = v + \mu_v$ & $\mu_v \sim P_v = 3e^{-7} $ \\ 
			$R_\textup{N}$ & $R_\textup{N} = R\textup{exp}(\mu_R) $& $\mu_R \sim P_R = 3e^{-8} $ \\ 
			$\Omega_\textup{N}$ & $\Omega_\textup{N} = \Omega + \mu_\Omega $ & $\mu_\Omega \sim P_\Omega = 3e^{-7} $ \\
			\hline
		\end{tabular}
	\end{center}
	\caption{Measurement noise level in power spectral density for the comparisons between LESO, FxTSDO, and FFTS-ESO}
	\label{table:ESO Comparison Noise Level}
\end{table}
\begin{figure}[htbp]
	\begin{minipage}[t]{0.5\linewidth}
		\centering
		\scriptsize
		\subfloat[Hover]{
			\includegraphics[width=\columnwidth]{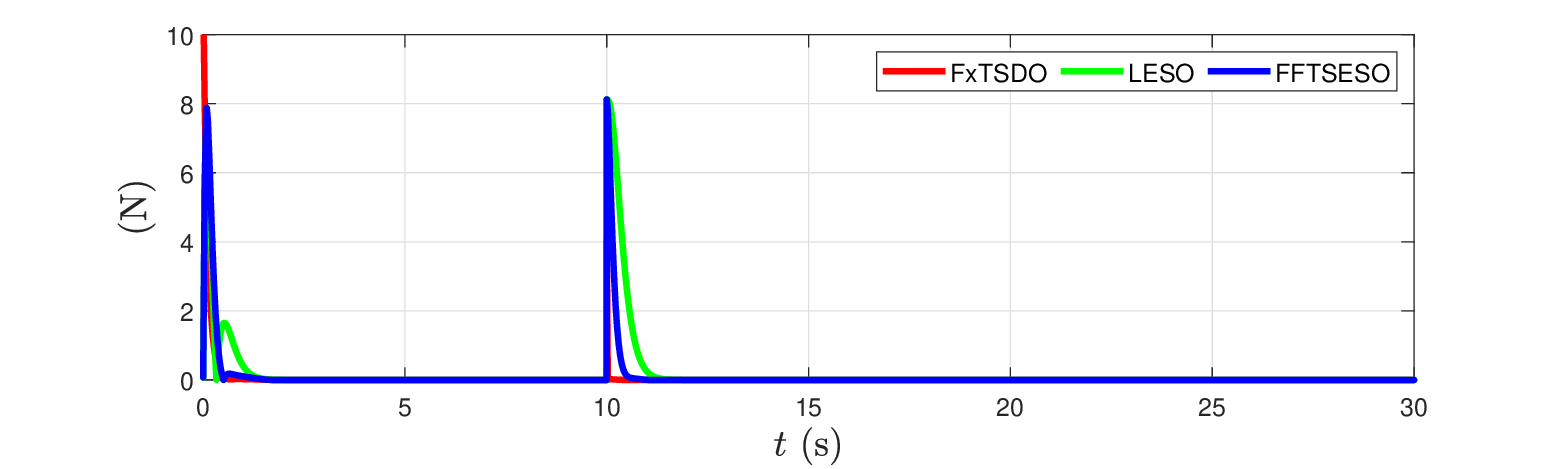}}\\
		\subfloat[Slow swing]{
			\includegraphics[width=\columnwidth]{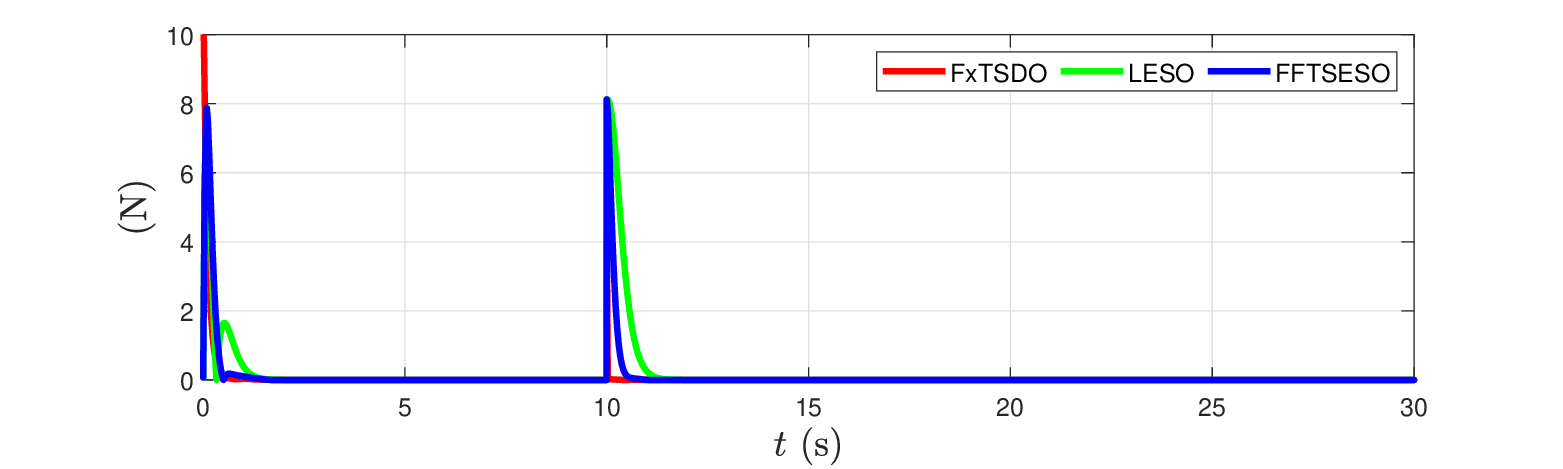}}\\
		\subfloat[Fast swing]{
			\includegraphics[width=\columnwidth]{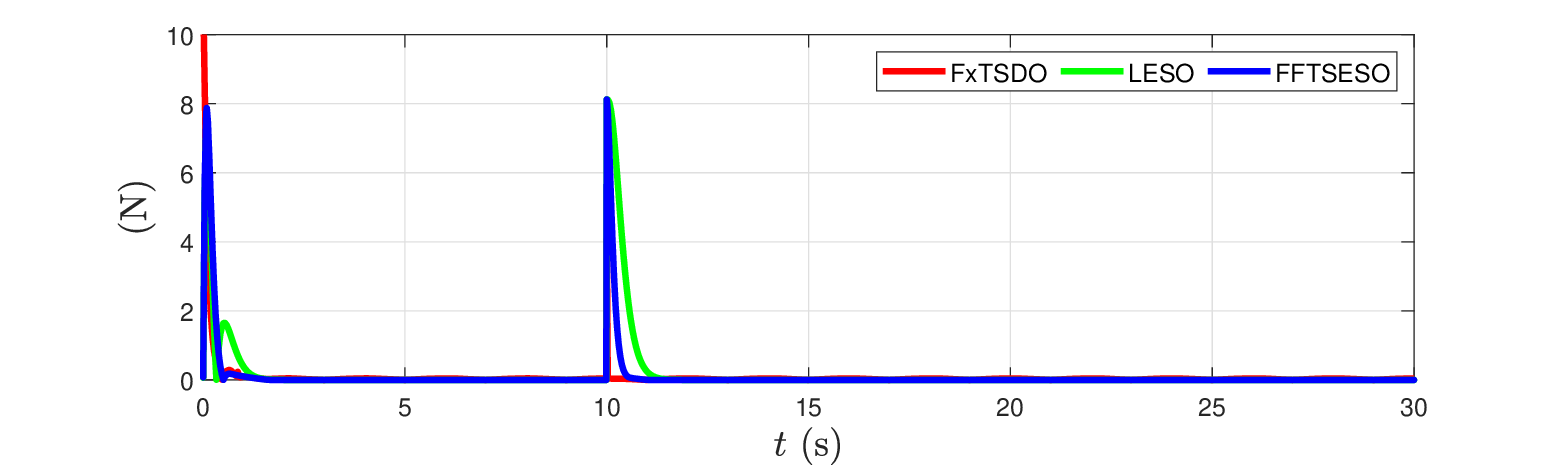}}\\
		\subfloat[High pitch]{
			\includegraphics[width=\columnwidth]{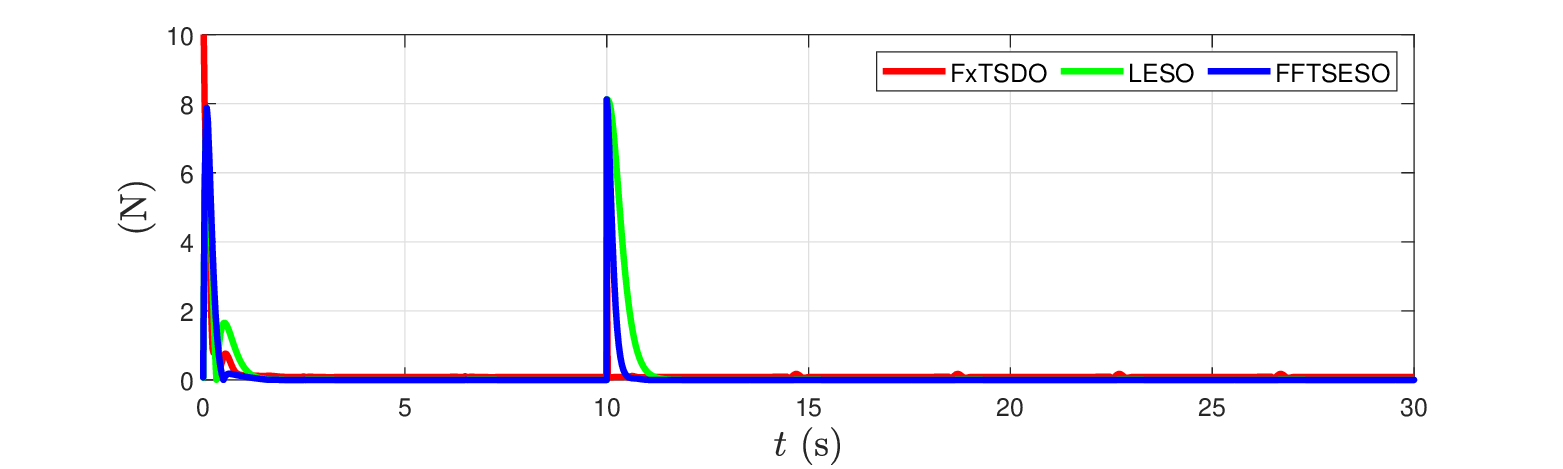}}
		\caption{Disturbance force estimation errors of the UAV from FxTSDO, LESO, and FFTS-ESO, in four different tracking control scenarios without measurement noise. }
		\label{fig:Disturbance Force DOESO Comparison Without Noise}
	\end{minipage}
	\begin{minipage}[t]{0.5\linewidth}
		\centering
		\scriptsize
		\subfloat[Hover]{
			\includegraphics[width=\columnwidth]{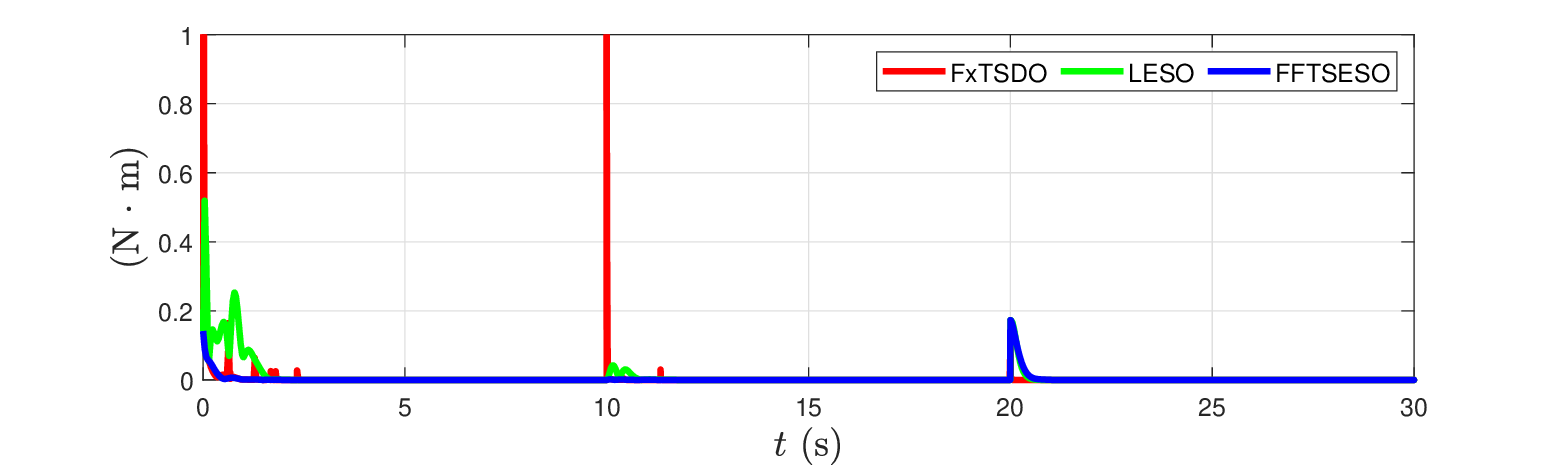}}\\
		\subfloat[Slow swing]{
			\includegraphics[width=\columnwidth]{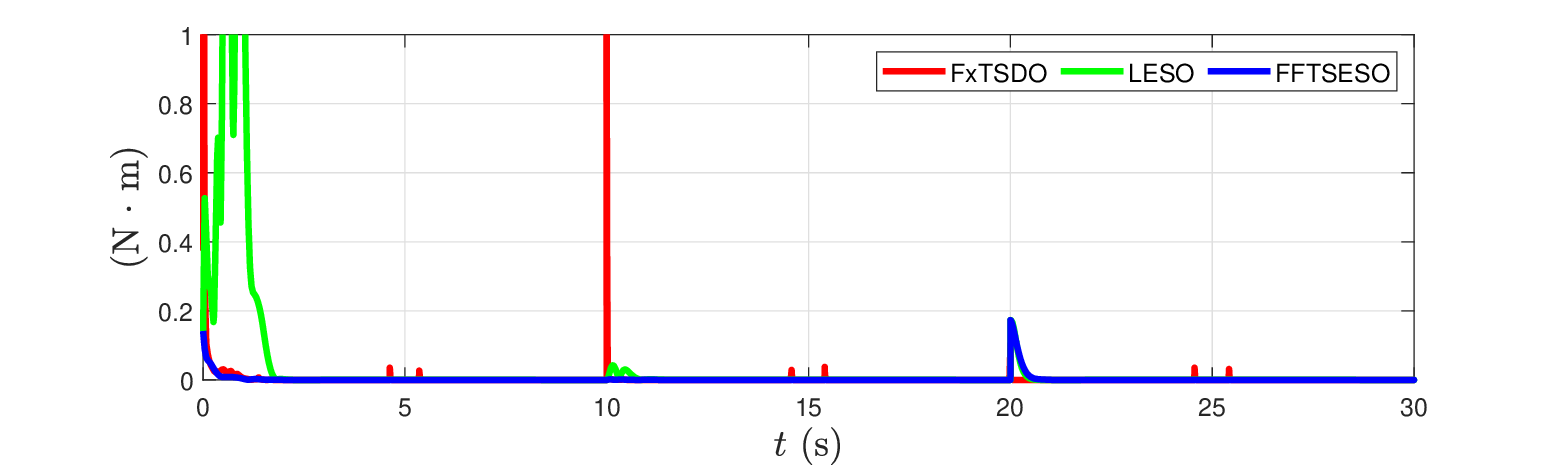}}\\
		\subfloat[Fast swing]{
			\includegraphics[width=\columnwidth]{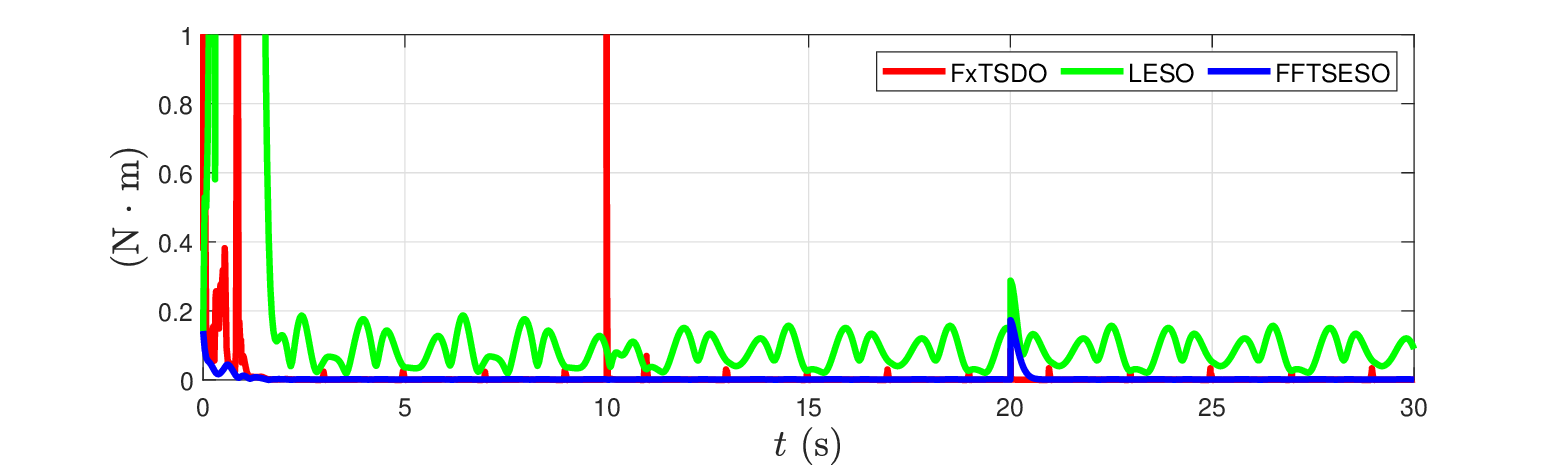}}\\
		\subfloat[High pitch]{
			\includegraphics[width=\columnwidth]{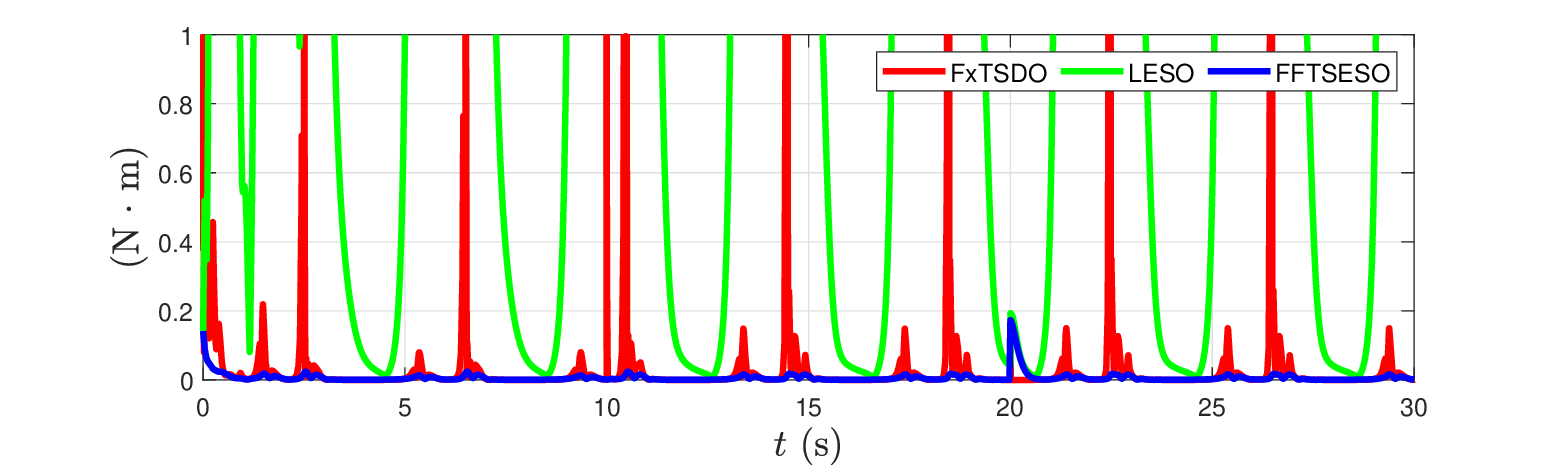}}
		\caption{Disturbance torque estimation errors of the UAV from FxTSDO, LESO, and FFTS-ESO, in four different tracking control scenarios without measurement noise.}
		\label{fig:Disturbance Torque DOESO Comparison Without Noise}
	\end{minipage}
\end{figure}

The four flight scenarios are the four desired trajectories listed in Table \ref{table:ESO Comparison Trajectory}. `Hovering' is the simplest flight scenario where the aircraft is ordered to hover at a fixed position during the simulation. `High Pitch' is the most complex flight scenario where the aircraft has to pitch up and 
track a circular trajectory. Since the norm of centripetal acceleration in the `High Pitch' scenario is more than a $g$, the aircraft has to flip over to track the desired trajectory. This desired trajectory with high centripetal acceleration forces the aircraft to go past the 90$^\circ$ pitch 
singularity of an Euler angle attitude representation. The measurement noise levels are as listed in Table \ref{table:ESO Comparison Noise Level} in terms of power spectral density (PSD). The disturbance force and 
torque in all of the four scenarios in this set of simulations are identical and they are the following step functions: 
\begin{align*}
	\varphi_{D} (t) &=
	\begin{cases}
		[5, \ 10, \ 0]\Tp \ \textup{N} & t< 10 \ \textup{s}\\
		[9, \ 15, \ 5]\Tp \ \textup{N} & t \geq 10 \ \textup{s} \\
	\end{cases}, 
	\tau_{D} (t) =
	\begin{cases}
		[-0.1, \ 0.1, \ 0.1]\Tp \ \textup{N}\cdot \textup{m} & t< 20 \ \textup{s} \\
		[0, \ 0, \ 0.2]\Tp \ \textup{N}\cdot \textup{m} & t \geq 20 \ \textup{s} \\
	\end{cases}
\end{align*} 

The parameters for FFTS-ESO in these simulations are $p=1.2, k_{t1}=3, k_{t2}=2,  k_{t3}=6, \kappa_t =0.8$, $k_{a1}=3, k_{a2}=2, k_{a3}=4, \kappa_a =0.6$. The gains for FxTSDO and LESO are as given by \cite{liu2022fixed} and by \cite{shao2018robust}. 
In the simulated flight, the initial states of the UAV for all four scenarios are: $R(0) = I, \ \Omega(0) = \left[0, \ 0, \ 0\right]\T \text{rad/s},  b(0) = \left[ 0.01, \ 0, \ 0 \right]\T \text{m}, \ v(0) = \left[5\pi, \ 0, \ 0\right]\T \text{m/s}.$ The initial conditions of the FxTSDO, LESO, and FFTS-ESO, are identical to the pose, velocities and disturbance of the UAV at the initial time in the simulation. 

\begin{figure}[htbp]
	\begin{minipage}[t]{0.5\linewidth}
		\centering
		\scriptsize
		\subfloat[Hover]{
			\includegraphics[width=\columnwidth]{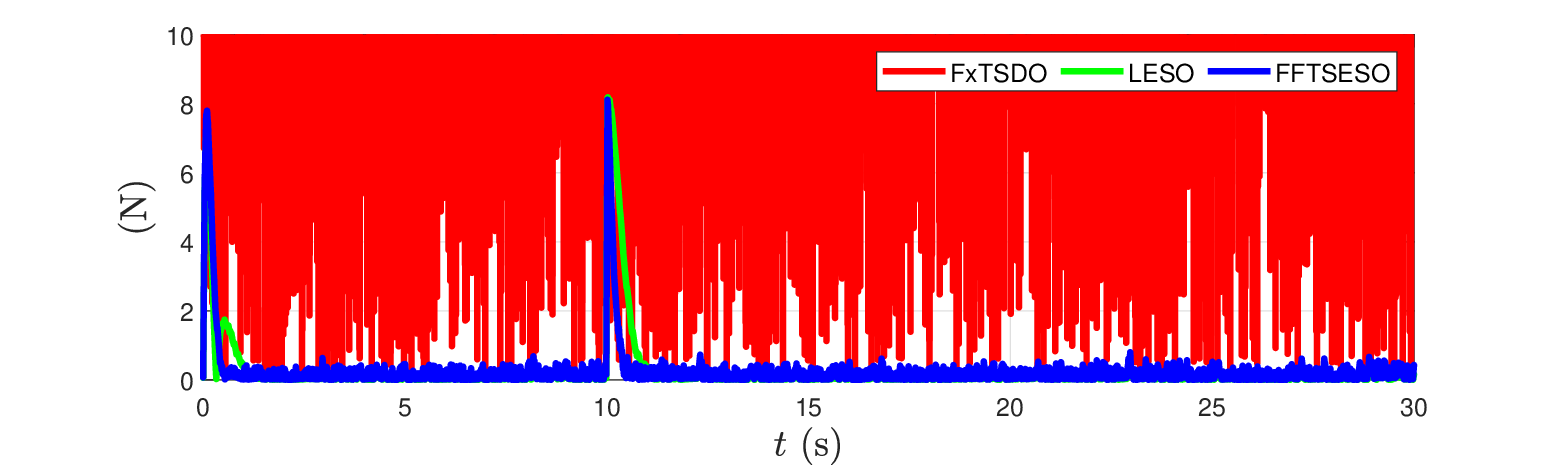}}\\
		\subfloat[Slow swing]{
			\includegraphics[width=\columnwidth]{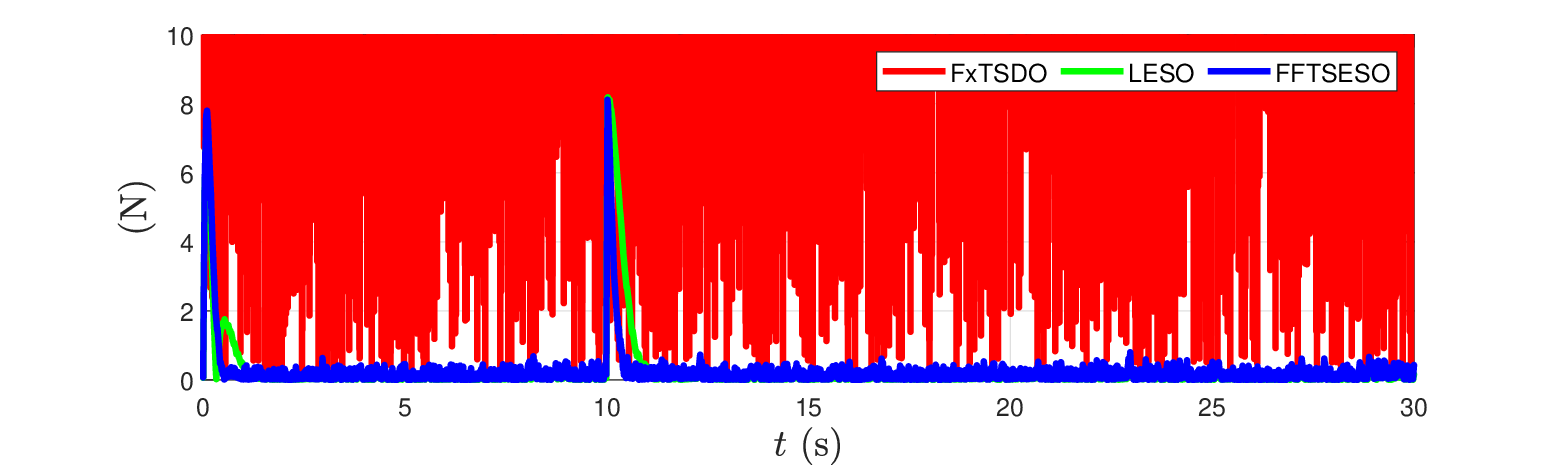}}\\
		\subfloat[Fast swing]{
			\includegraphics[width=\columnwidth]{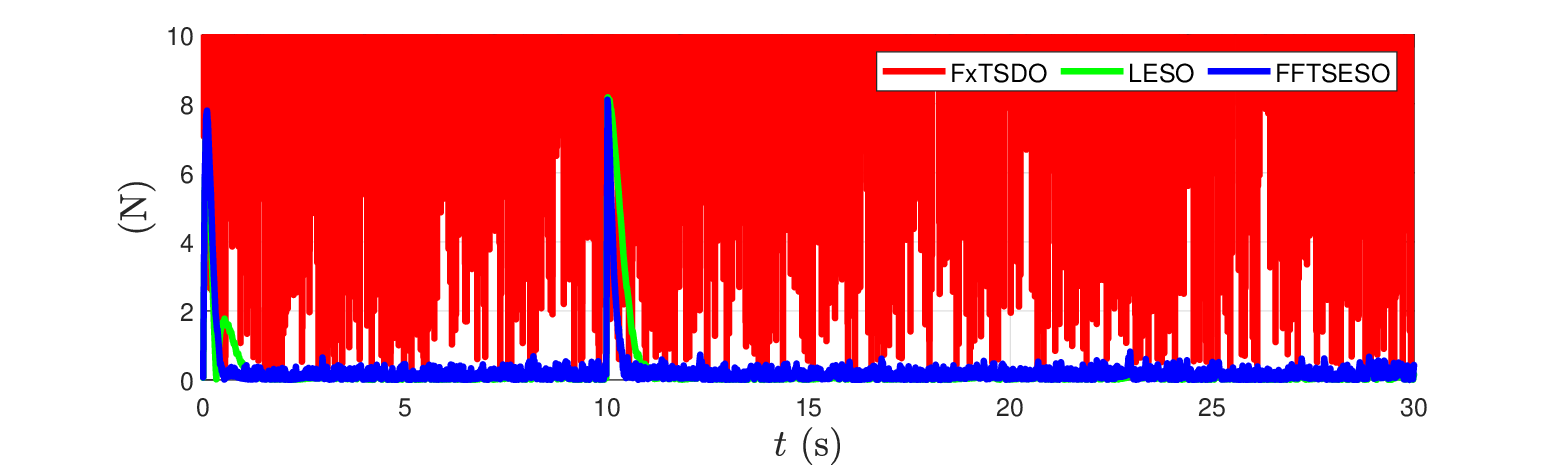}}\\
		\subfloat[High pitch]{
			\includegraphics[width=\columnwidth]{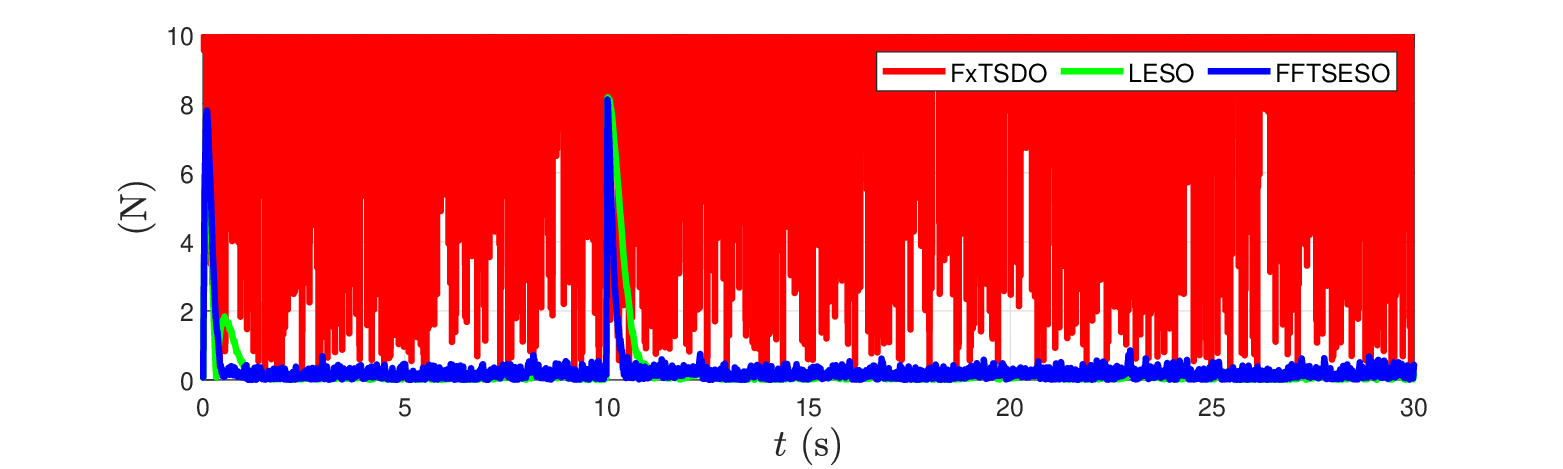}}
		\caption{Disturbance force estimation errors of the UAV from FxTSDO, LESO, and FFTS-ESO, in four different tracking control scenarios with measurement noise.}
		\label{fig:Disturbance Force DOESO Comparison With Noise}
	\end{minipage}
	\begin{minipage}[t]{0.5\linewidth}
		\centering
		\scriptsize
		\subfloat[Hover]{
			\includegraphics[width=\columnwidth]{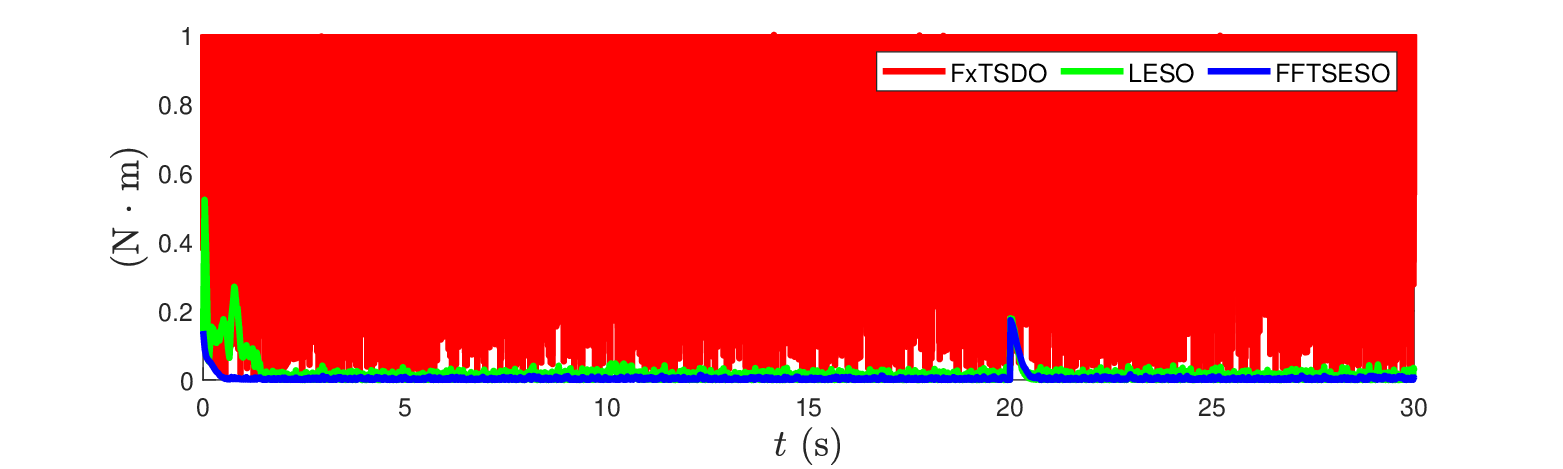}}\\
		\subfloat[Slow swing]{
			\includegraphics[width=\columnwidth]{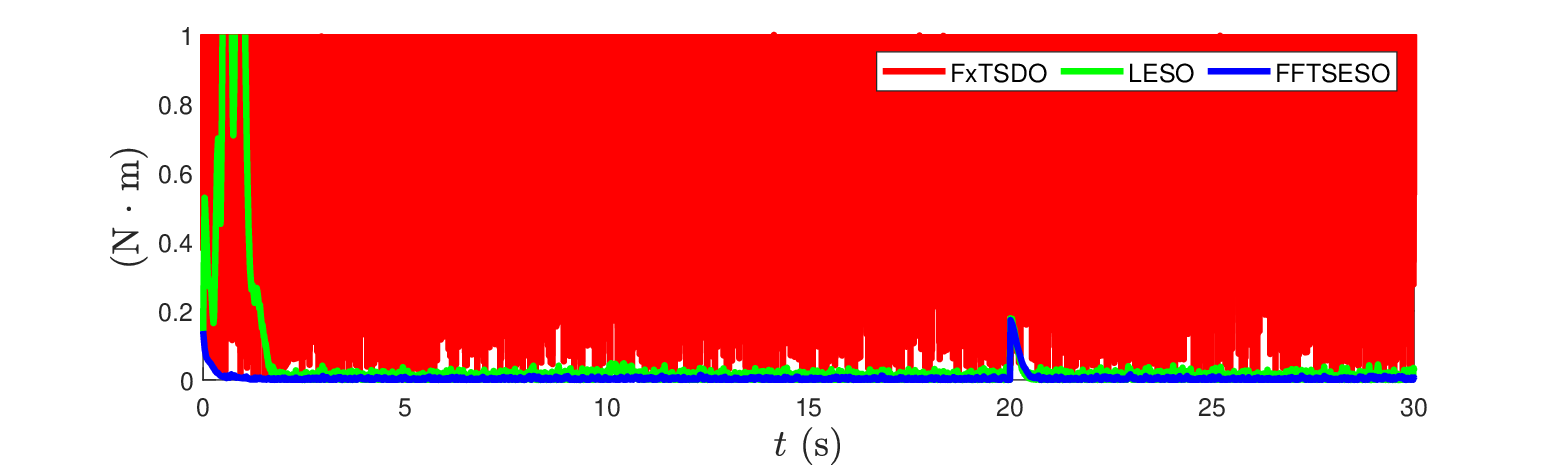}}\\
		\subfloat[Fast swing]{
			\includegraphics[width=\columnwidth]{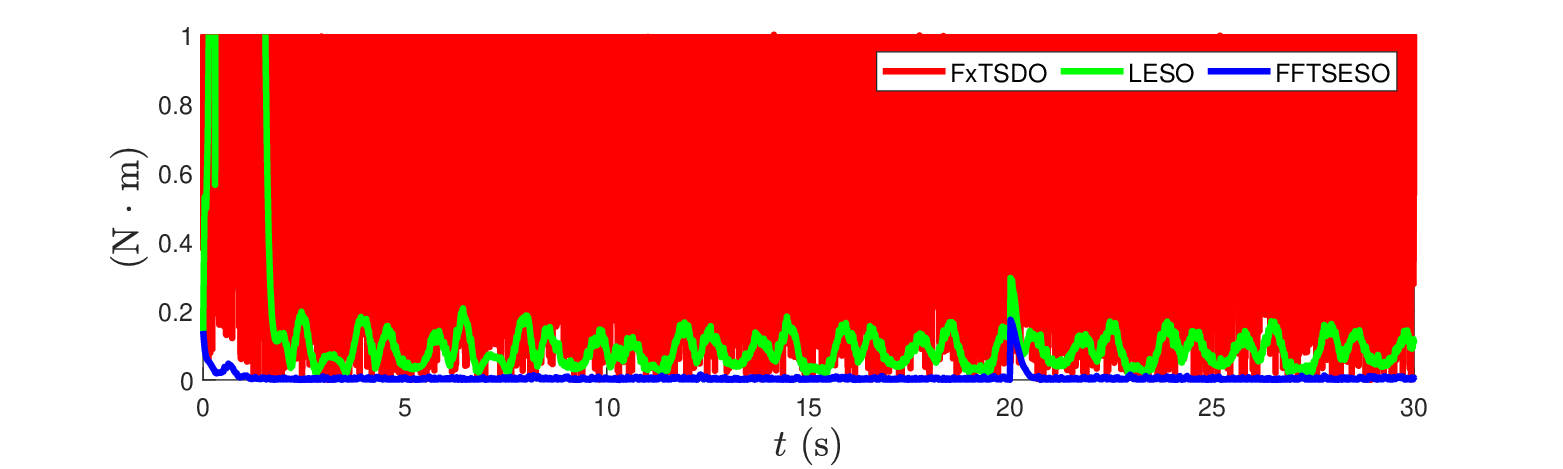}}\\
		\subfloat[High pitch]{
			\includegraphics[width=\columnwidth]{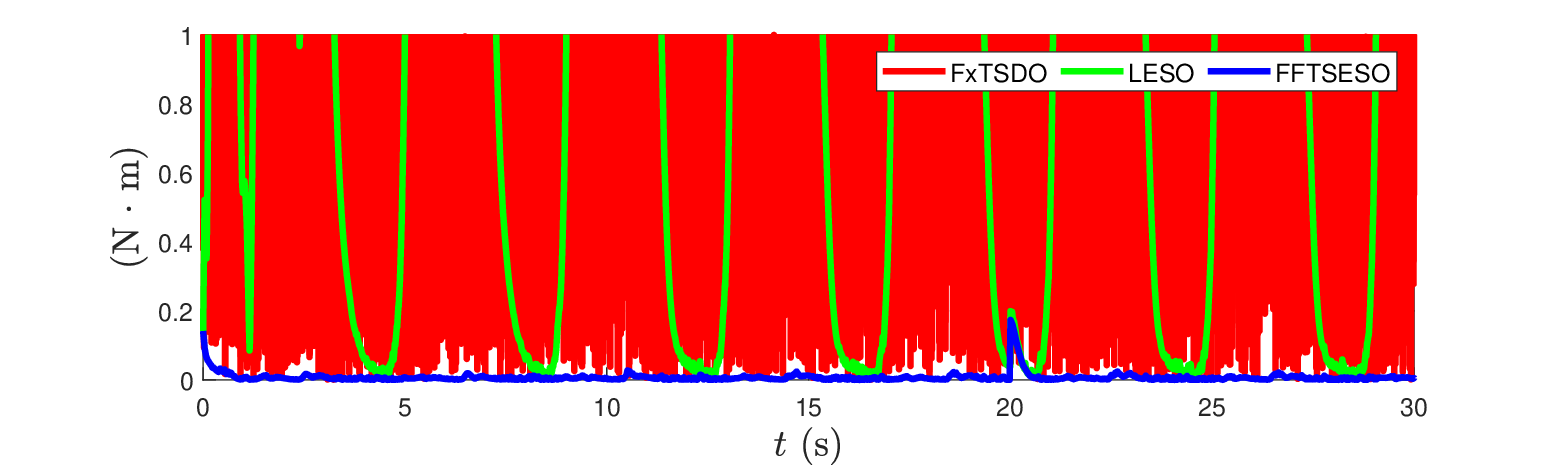}}
		\caption{Disturbance torque estimation errors of the UAV from FxTSDO, LESO, and FFTS-ESO, in four different tracking control scenarios with measurement noise.}
		\label{fig:Disturbance Torque DOESO Comparison With Noise}
	\end{minipage}
\end{figure}

We present the simulation results in four sets of figures. Figure \ref{fig:Disturbance Force DOESO Comparison Without Noise} and \ref{fig:Disturbance Torque DOESO Comparison Without Noise} present the disturbance force and torque estimation errors respectively, from FxTSDO, LESO and FFTS-ESO in the flight scenarios described in Table \ref{table:ESO Comparison Trajectory} with noise-free measurements. Figure \ref{fig:Disturbance Force DOESO Comparison With Noise} and \ref{fig:Disturbance Torque DOESO Comparison With Noise} present the disturbance estimation errors from these schemes for the flight trajectories in Table \ref{table:ESO Comparison Trajectory}, in the presence of  measurement noise levels as described in Table \ref{table:ESO Comparison Noise Level}.

Figure \ref{fig:Disturbance Force DOESO Comparison Without Noise} shows the disturbance force estimation errors from the three schemes with noise-free measurements. Although the disturbance force estimation error from FxTSDO shows significant initial transient, the results from Figure \ref{fig:Disturbance Force DOESO Comparison Without Noise} indicates that with noise-free measurement, the disturbance force estimations from these three schemes converge to the origin in all four flight scenarios. The transients at $t=15$ s are from the step-function disturbance force $\varphi_D$, whose step time is $t=15$ s.
Figure \ref{fig:Disturbance Torque DOESO Comparison Without Noise} shows the disturbance torque estimation errors from the three schemes with noise-free measurement. In Figure \ref{fig:Disturbance Torque DOESO Comparison Without Noise}, we observe that when $t=10 \ \text{s}$, high transients appears in the disturbance torque estimation error from FxTSDO.

Despite the initial transients, the disturbance torque estimation errors from all three schemes converge to the origin in 'Hovering' and 'Slow swing' scenarios. However, in 'Fast swing' and 'High pitch' scenarios, the disturbance torque estimation errors from LESO and FxTSDO diverge. As is stated in Section \ref{sec:Intro}, since the LESO uses Euler-angle to represent attitude for disturbance torque estimation, it experiences a singularity in attitude representation when the UAV tracks the 'Fast swing' and 'High Pitch' trajectories. Thus, in these two scenarios, the singularity in the attitude representation destabilizes the disturbance torque estimation error of LESO. 

Figure. \ref{fig:Disturbance Force DOESO Comparison With Noise} and \ref{fig:Disturbance Torque DOESO Comparison With Noise} present the disturbance force and disturbance torque estimation errors respectively, from the three schemes with identical noisy measurements as given in Table \ref{table:ESO Comparison Noise Level}. As is stated in Remark \ref{rem:Compare}, we observe that with measurement noise, FxTSDO is not capable of providing any meaningful disturbance estimation.  In 'Fast swing' and 'High pitch' scenarios, the disturbance torque estimation errors from LESO diverge from the origin. 

To summarize, figures \ref{fig:Disturbance Force DOESO Comparison Without Noise}, \ref{fig:Disturbance Torque DOESO Comparison Without Noise}, \ref{fig:Disturbance Force DOESO Comparison With Noise}, and \ref{fig:Disturbance Torque DOESO Comparison With Noise} show that the FFTS-ESO has satisfactory disturbance estimation performance and outperforms the LESO and FxTSDO when the UAV experiences large pose changes and has noisy measurements.

\section{Flight Experiments}\label{sec:Experiment}
In this section, the proposed FFTS-ESO is validated through flight experiments. Its hardware and software are custom-designed and developed based on the open-source autopilot PX4 by \cite{meier2015px4}. To demonstrate the capability of estimating and rejecting 
the disturbances, flight experiments are conducted under wind disturbances generated by a fan arrary wind tunnel (FAWT) from the Switzerland-based company WindShape. We first describe the hardware and software 
configurations of the UAV and the setup of the experiment. Afterwards, we present our experimental results including the characteristics of the wind disturbances and 
the control performance of the UAV when exposed to disturbances generated by the FAWT.

\subsection{Hardware configuration}
\begin{figure}[htbp]
	\begin{minipage}[t]{0.48\linewidth}
		\centering
		\scriptsize
		\includegraphics[width=0.83\columnwidth]{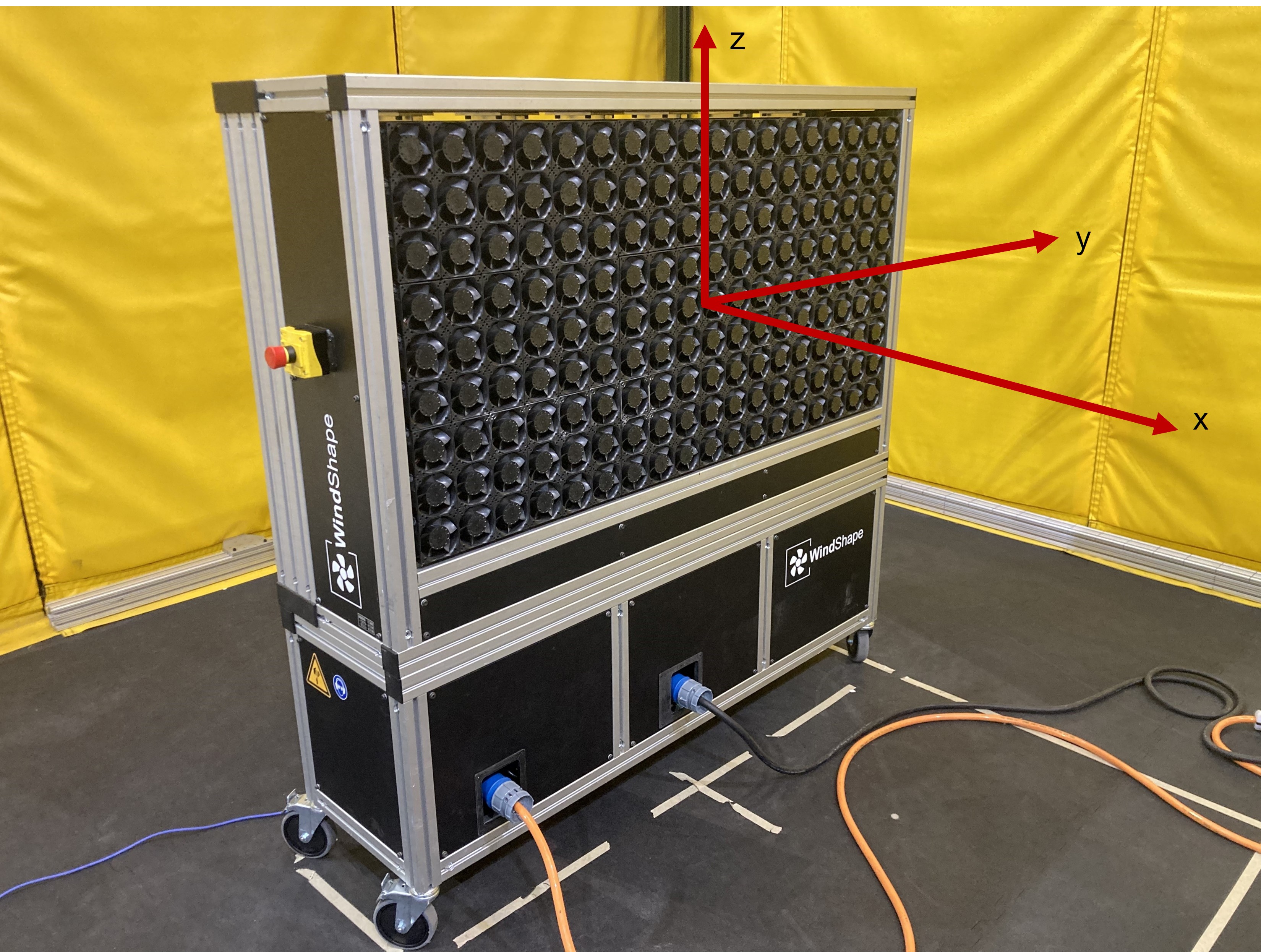}
		\caption{The FAWT and its body-fixed coordinate frame.}
		\label{fig:FAWT}
	\end{minipage}
	\begin{minipage}[t]{0.48\linewidth}
		\centering
		\scriptsize
		\includegraphics[width=0.83\columnwidth]{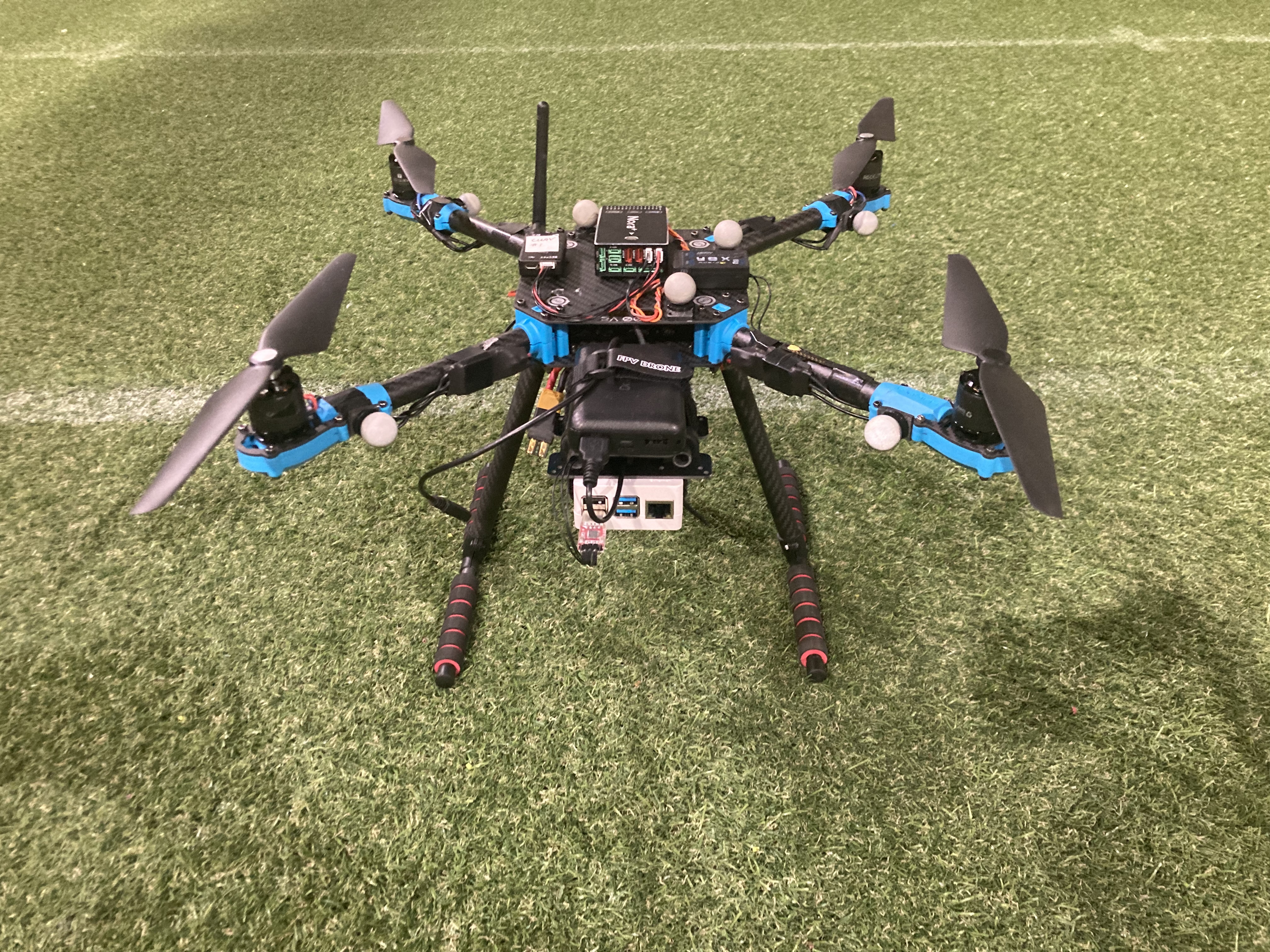}
		\caption{Quadrotor UAV developed in the autonomous unmanned system lab (AUSL) at Syracuse University.}
		\label{fig:Aircraft}
	\end{minipage}
\end{figure}

The multi-rotor UAV is shown in Figure \ref{fig:Aircraft}. It has four brushless direct current electrical motors (T-Motor Air 2216 880-KV) paired with
$10'' \times 4.5''$ carbon fiber propellers. To control the rotational speed of the motors, each is connected to an electronic speed control (T-Motor Air 20A), which receives commands from a PixHawk flight control unit (FCU, CUAV Nora plus) with redundant inertial measurement units (IMU). Flight control and state estimation are conducted by the FCU. The pose of the vehicle is measured by an optical motion capture system (VICON), and sent to a companion computer (Raspberry Pi 4) through Wi-Fi network, and then passed on to the FCU through a telemetry port.
To generate wind fields with various turbulence characteristics for the flight experiments, we leveraged the Fan Array Wind Tunnel (FAWT) in Figure. \ref{fig:FAWT}. The wind tunnel is a
1.44 m by 0.72 m array of 162 independently controllable fans capable of generating wind speeds up to 12m/s. The distributed fans are controlled in real time by a Python program. 
In our flight experiments, the FAWT was run at the steady uniform flow mode, with each individual fan running at identical and constant duty. 

\subsection{Software configuration}
The flight control software is developed from the open-source autopilot software PX4 v1.13.2. According to \cite{meier2015px4}, the system architecture of PX4 is centered around a publish-subscribe object request broker on top of a POSIX application
programming interface. This programming interface has different modules for data logging, communication, estimation, and control. The FFTS-ESO is implemented onto the module \texttt{mc\_pos\_control} and \texttt{mc\_rate\_control} for translational 
and rotational motions, respectively. The feedback of disturbance estimates from the FFTS-ESO is applied to the control law as an additional term, so that the original control architecture is modified with this feedforward disturbance rejection term. We 
introduce Boolean parameters to switch the disturbance rejection conveniently. 

In the experiment, the rest of the autopilot (PX4 v1.13.2) is kept unchanged, to have a fair comparison of the flight control performance between the original PX4 autopilot, and the one with disturbance rejection from FFTS-ESO. The flight control parameters of the autopilot are as described in the multi-rotor frame S500 in \href{https://github.com/PX4/PX4-Autopilot/blob/98d893503495f7c28856bccf830082451b20265d/ROMFS/px4fmu_common/init.d/airframes/4014_s500#L4}{https://github.com/PX4/PX4-Autopilot/}. A Robot Operating System (ROS) interface program is developed for the companion computer that transmits commands and pose to the vehicle. The flight data are saved in the memory card inside the FCU in the form of .ulg file for post-processing. We use the MAVLINK telecommunication protocol for communication between the FCU, companion computer, and ground control station. 

The FFTS-ESO parameters are selected as: $p=1.2, k_{t1}=6, k_{t2}=3,  k_{t3}=1, \kappa_t =0.6$, $k_{a1}=8, k_{a2}=4, k_{a3}=2, \kappa_a =0.6$. The empirically known mass and inertia of the vehicle as given to the FFTS-ESO are: $m=1$ kg and $  J=\textup{diag}([0.03,0.03,0.06]) \ \textup{kg}\cdot \textup{m}^2$.
We link the source code of the customized PX4 with FFTS-ESO on Github.\footnote{Github link: \url{https://github.com/nswang1994/GeometricPX4/tree/Geometric-FFTS-ESO}} 

\subsection{Experiment procedure}
\begin{figure}[htbp]
	\centering
	\scriptsize
	\includegraphics[width=0.6\columnwidth]{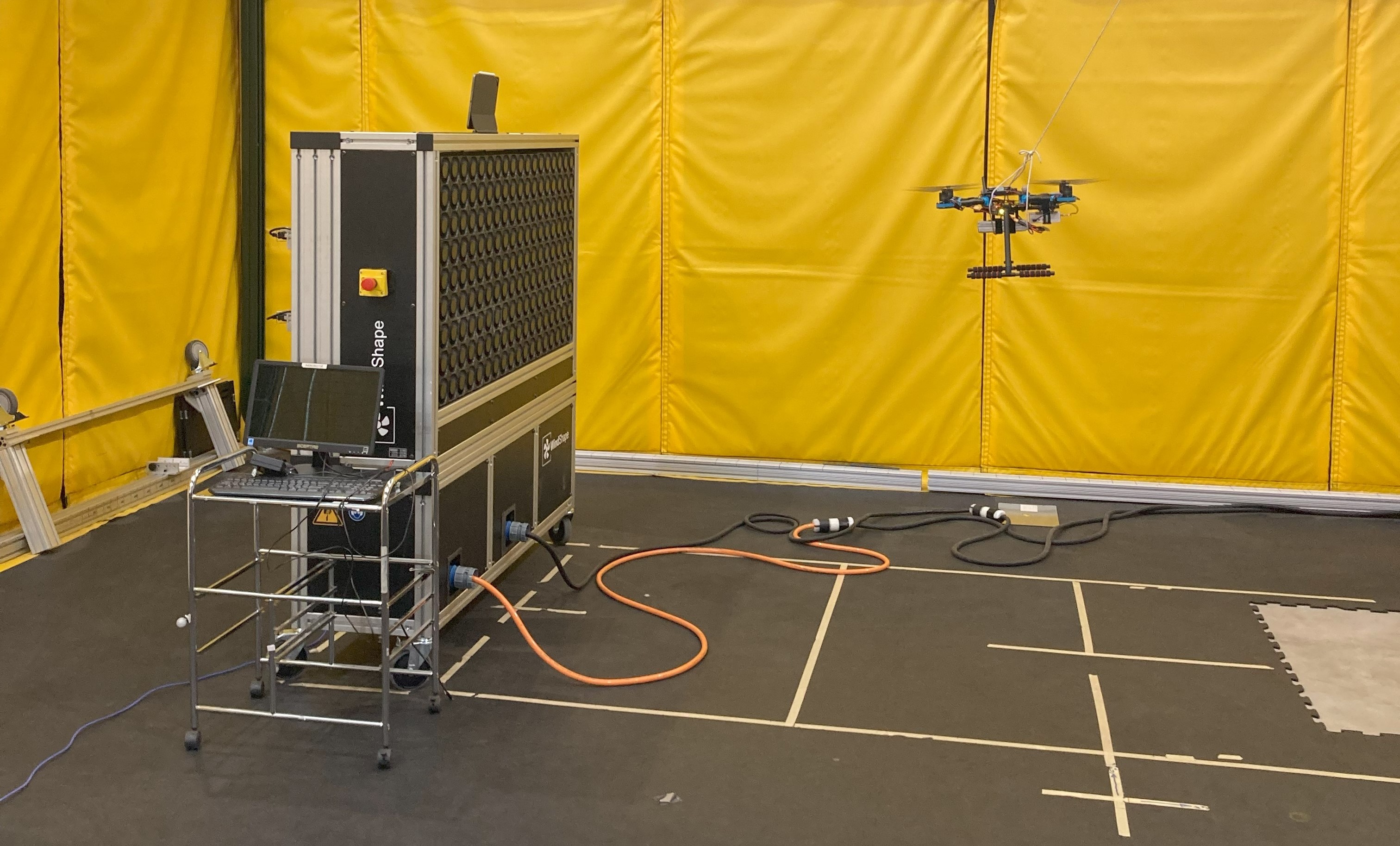}
	\caption{Setup for our flight control experiments.}
	\label{fig:Flight}
\end{figure}

The flight experiment setup is shown in Figure \ref{fig:Flight}. We define the FAWT coordinate frame as shown in Figure \ref{fig:FAWT}, with $x$ as the stream-wise direction, $y$ as the span-wise direction, and $z$ as the vertically up direction. The origin is at the geometric center of the fan array. 
We operate the FAWT at a steady uniform flow mode at 30\%-70\% of its maximum duty to measure the wind velocity of the wind field. We conduct the wind velocity measurements with a hotwire anemometer facing in the $x$ direction at $x=1.2\text{m}, y=0\text{m}, z=0\text{m}$ in the FAWT coordinate system. 

As shown in Figure \ref{fig:Flight}, the vehicle is commanded to hover in the front of the FAWT, at $x=1.5\text{m}, y=0\text{m},  z=0\text{m}$ in the FAWT frame. This 
hovering position is at the center point of the test section, so that we can maximally avoid the boundary layer around the section border, where higher turbulence intensity and flow uncertainty occur. The time for hovering flight is set to 210 s. During this period, we turn on the FAWT for 150 s to disturb the 
vehicle with turbulent flows with statistically constant characteristics. The pose of the vehicle during flight is recorded in the log file for evaluation. 

\subsection{Results: turbulent flow measurement}
The results for hot-wire measurements are covered in this subsection. For brevity, we omit the details of the measurement procedure, which are accessible in the dissertation by \cite{wang2023geometric}. 
According to the Reynolds decomposition given in \cite{pope2000turbulent}, we decompose the measured wind velocity $u$ into the sum of a time-averaged velocity $\overline{u}$ and fluctuating velocity $\widetilde{u}$, such 
that $u = \overline{u} + \widetilde{u}$, where $\overline{(\cdot)} $ stands for the time-averaged quantity. Based on hotwire measurements, we characterize the FAWT wind field with the following quantities: the time-averaged velocity $\overline{u}$, the variance of the fluctuating velocity $\overline{\widetilde{u}^2}$, and the turbulence intensity (TI), which is defined by $\sqrt{\overline{\widetilde{u}^2}}/\overline{u}$. We itemize the results as follows:
\begin{itemize}
	\item $\overline{u_{30}} =5.472 \ \text{m}/\text{s}$; $\overline{\widetilde{u_{30}}^2} = 0.061 \ \text{m}^2/\text{s}^2$; $\text{TI}_{30} = 0.0451$;
	\item $\overline{u_{40}} =6.876 \ \text{m}/\text{s}$; $\overline{\widetilde{u_{40}}^2} = 0.082 \ \text{m}^2/\text{s}^2$; $\text{TI}_{40} = 0.0417$;
	\item $\overline{u_{50}} =8.213 \ \text{m}/\text{s}$; $\overline{\widetilde{u_{50}}^2} = 0.116 \ \text{m}^2/\text{s}^2$; $\text{TI}_{50} = 0.0415$;
	\item $\overline{u_{60}} =9.590 \ \text{m}/\text{s}$; $\overline{\widetilde{u_{60}}^2} = 0.168 \ \text{m}^2/\text{s}^2$; $\text{TI}_{60} = 0.0427$;
	\item $\overline{u_{70}} =10.920 \ \text{m}/\text{s}$; $\overline{\widetilde{u_{70}}^2} = 0.237 \ \text{m}^2/\text{s}^2$; $\text{TI}_{70} = 0.0446$.
\end{itemize}
To summarize, for the FAWT in the setup described earlier in this subsection, we observe that $\overline{u}$ ranges from $5.472$ m/s to $10.920$ m/s, and TI is around $0.043$. Moreover, we observe that 
$\overline{\widetilde{u}^2}$ has positive correlation with $\overline{u}$. We assume that higher $\overline{\widetilde{u}^2}$ brings higher turbulence energy, which causes higher disturbance inputs to the aircraft flying within the wind field.

\subsection{Results: flight experiment}
\begin{figure}[htbp]
	\begin{minipage}[t]{0.5\linewidth}
		\centering
		\scriptsize
		\subfloat[Duty 30\%]{
			\includegraphics[width=\columnwidth]{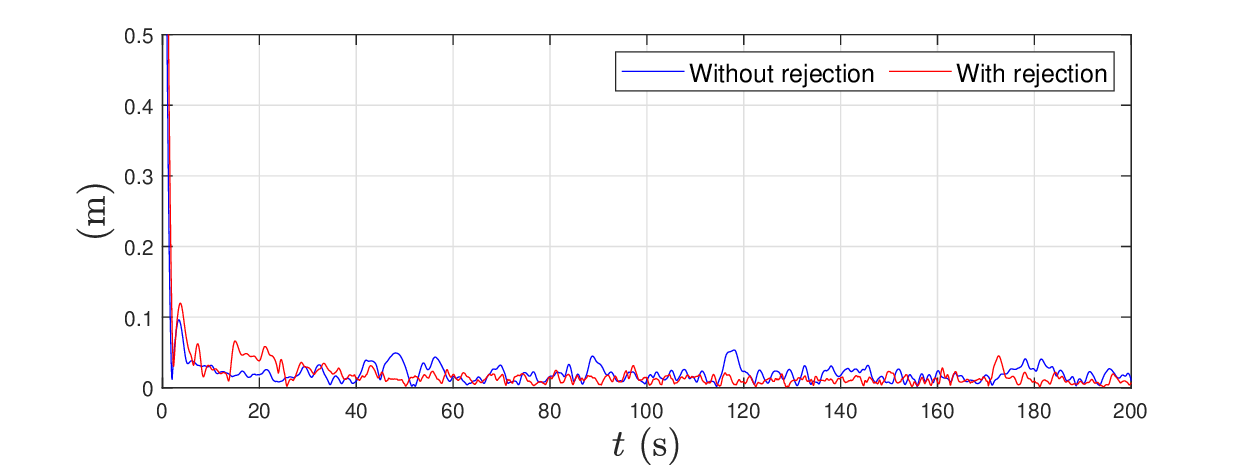}}\\
		\subfloat[Duty 40\%]{
			\includegraphics[width=\columnwidth]{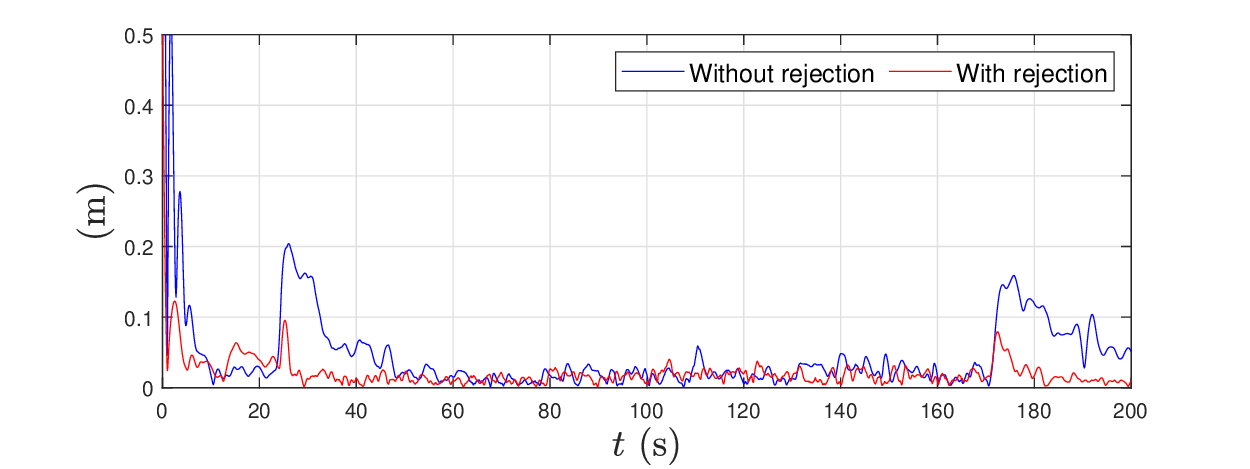}}\\
		\subfloat[Duty 50\%]{
			\includegraphics[width=\columnwidth]{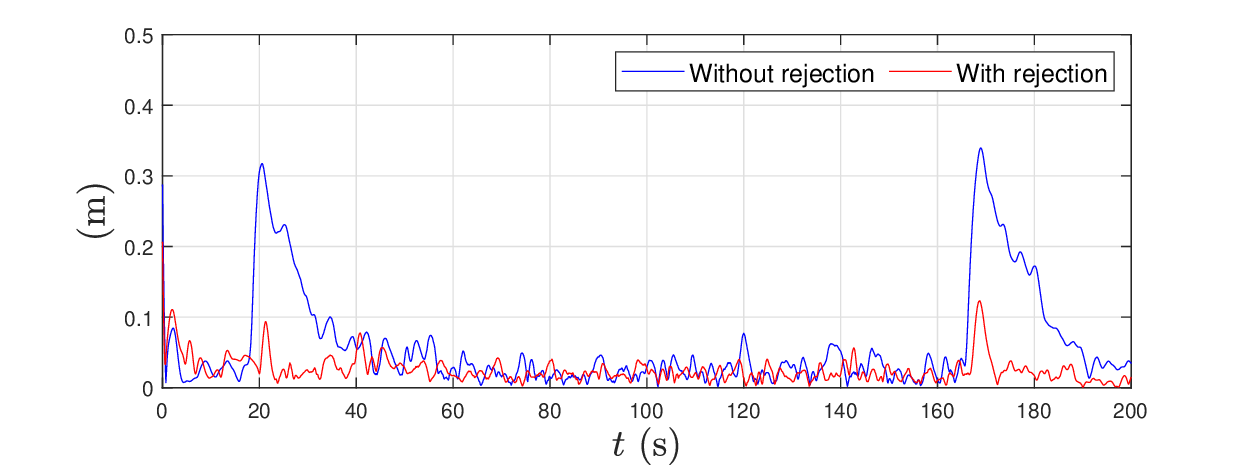}}\\
		\subfloat[Duty 60\%]{
			\includegraphics[width=\columnwidth]{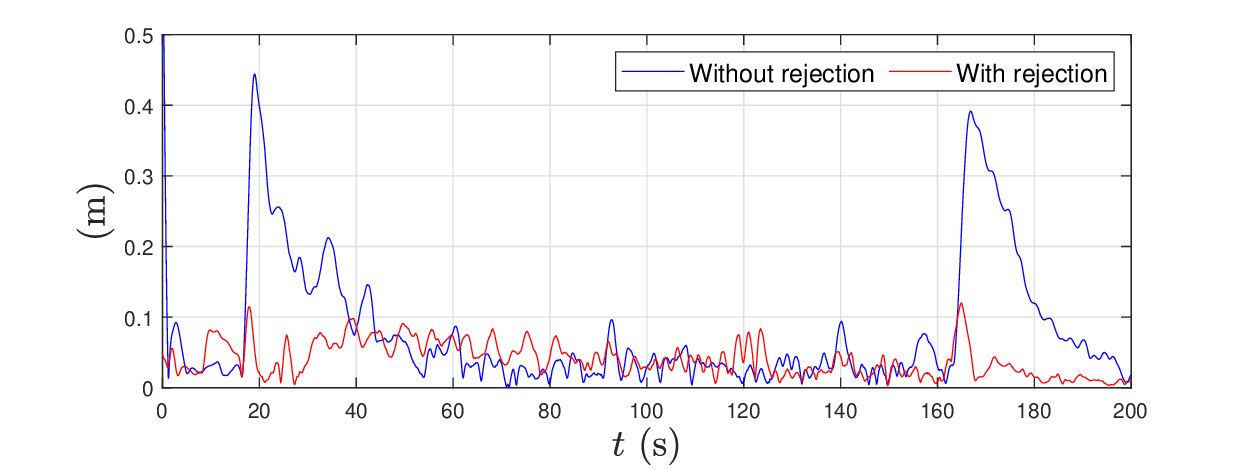}}\\
		\subfloat[Duty 70\%]{
			\includegraphics[width=\columnwidth]{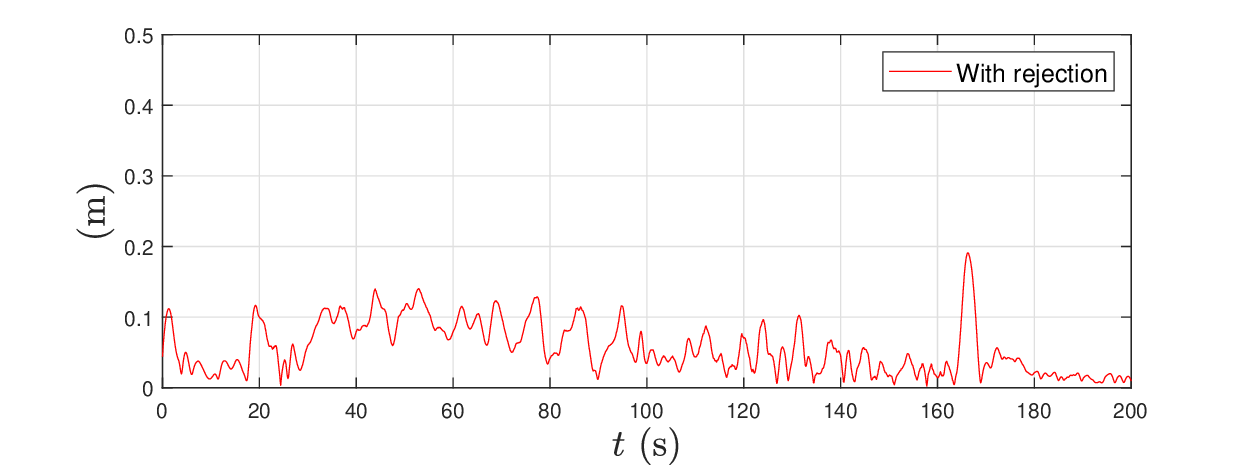}}
		\caption{Position tracking error}
		\label{fig:Position Tracking Error}
	\end{minipage}
	\begin{minipage}[t]{0.5\linewidth}
		\centering
		\scriptsize
		\subfloat[Duty 30\%]{
			\includegraphics[width=\columnwidth]{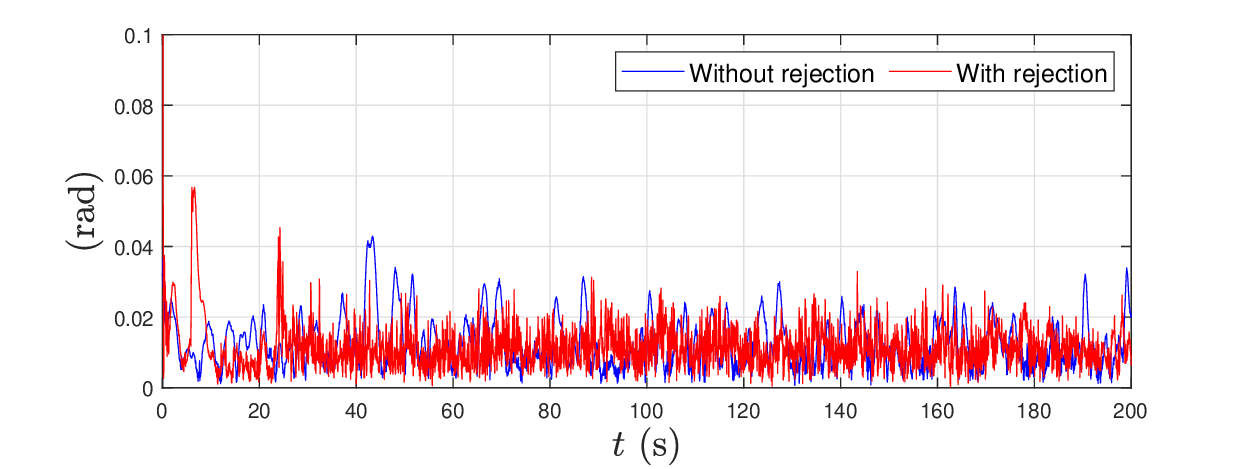}}\\
		\subfloat[Duty 40\%]{
			\includegraphics[width=\columnwidth]{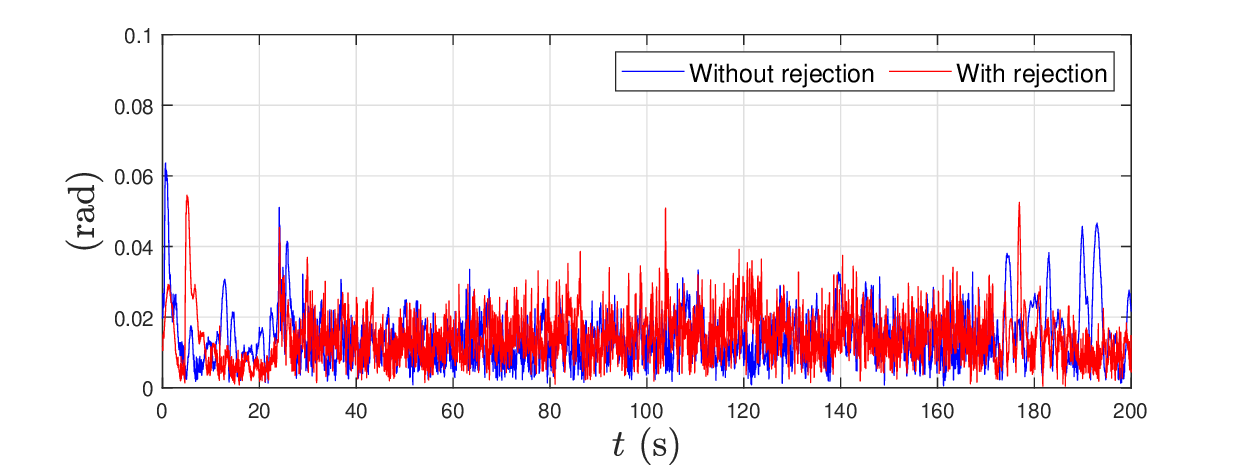}}\\
		\subfloat[Duty 50\%]{
			\includegraphics[width=\columnwidth]{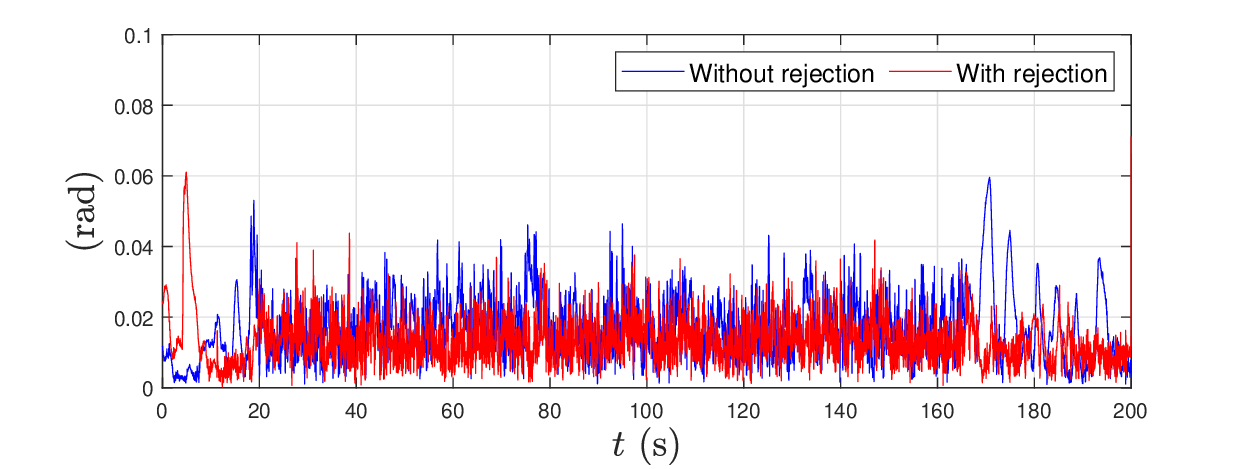}}\\
		\subfloat[Duty 60\%]{
			\includegraphics[width=\columnwidth]{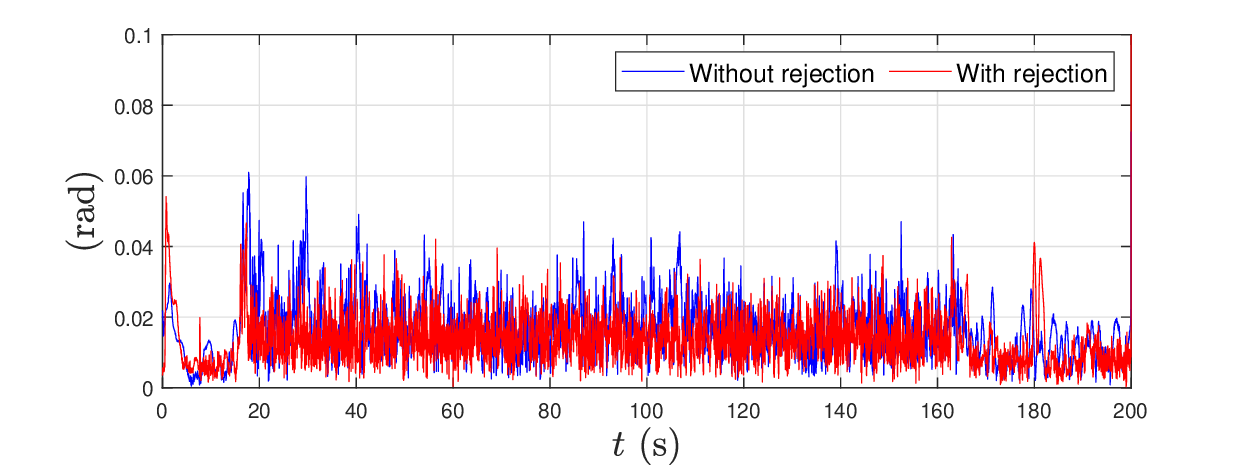}}\\
		\subfloat[Duty 70\%]{
			\includegraphics[width=\columnwidth]{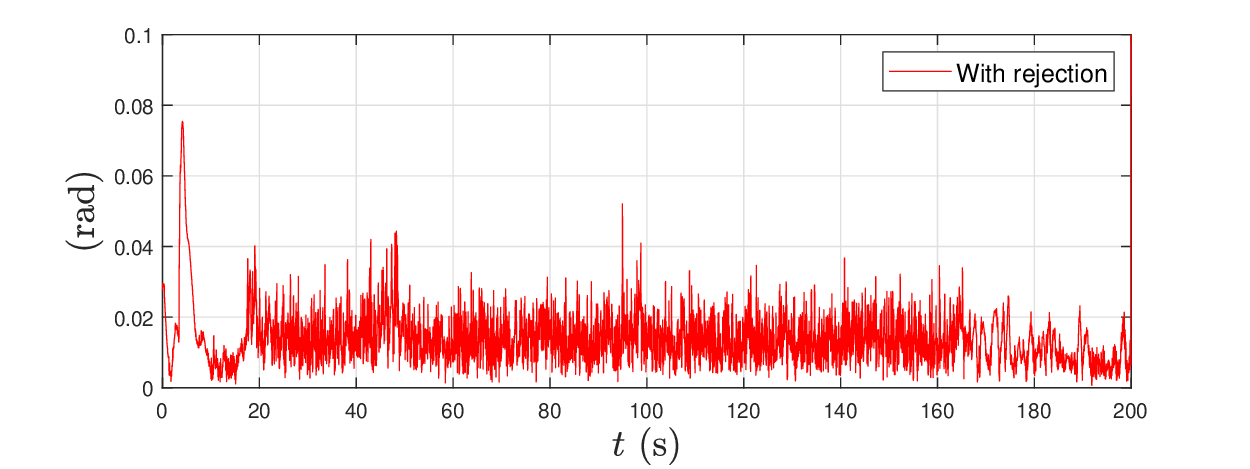}}
		\caption{Attitude tracking error}
		\label{fig:Attitude Tracking Error}
	\end{minipage}
\end{figure}

The results of the flight experiment are covered in this subsection. Figure \ref{fig:Position Tracking Error} and Figure \ref{fig:Attitude Tracking Error} illustrate the position and attitude tracking errors during hovering. The position tracking error is quantified as the Euclidean norm, while the attitude tracking error is quantified as the principal angle. To highlight the flight control performance under disturbances, we omit the pose data during take-off and landing in the presented result. The curves for the control scheme
without disturbance rejection are plotted in blue, and those with disturbance rejection are in red. The time-averaged position and attitude tracking errors are listed in Table \ref{table:Time-Averaged Tracking Error} .

\begin{table}[htbp]
	\begin{center}
		\begin{tabular}{ | c| c| c| c| c| } 
			\hline
			\ & \multicolumn{2}{c|}{Position tracking error (m)}  & \multicolumn{2}{c|}{Attitude tracking error (rad)} \\
			\hline
			\ & PX4 Stack & PX4+FFTS-ESO &PX4 Stack & PX4+FFTS-ESO  \\ 
			\hline
			30\% &  0.0251  & 0.0236  &  0.0125  & 0.0116 \\ 
			40\% &  0.0468  & 0.0211  &  0.0140  & 0.0141 \\ 
			50\% &  0.0589  & 0.0254  &  0.0166  & 0.0132 \\
			60\% &  0.0792  & 0.0400  &  0.0164  & 0.0134  \\
			70\% & {\color{red}Failed!} & 0.0557  &  {\color{red}Failed!} & 0.0139 \\
			\hline
		\end{tabular}
	\end{center}
	\caption{Time-averaged tracking error}
	\label{table:Time-Averaged Tracking Error}
\end{table}

Figures \ref{fig:Position Tracking Error} and \ref{fig:Attitude Tracking Error} show that both position and attitude tracking errors have high transient at around 20s and 180s when the disturbances from FAWT kick in and fade off, respectively. In Figure \ref{fig:Attitude Tracking Error}, for the attitude tracking error of the control scheme with disturbance rejection,  we observe extra transient at around 0s-10s, when the disturbance rejection kick-in. In Figure \ref{fig:Position Tracking Error}, we observe that when the FAWT operates at 40\%-60\% of its maximum duty, the position tracking error of the control scheme with disturbance rejection outperforms the one without rejection.  When the FAWT operates at 30\% of its maximum duty, the difference between the two control schemes is not evident in Figure \ref{fig:Position Tracking Error}. However, in terms of the time-averaged position tracking errors in Table \ref{table:Time-Averaged Tracking Error}, we can still observe that the scheme with disturbance rejection outperforms the one without rejection when the FAWT operates at 30\%-60\% of its maximum duty. When the FAWT operates at 70\% of its maximum duty, the control scheme without disturbance rejection mechanism fails to hover constantly, while the one with rejection succeeds. 

\section{Conclusion}\label{sec:Conclusion}
In this article, FFTS-ESO for disturbance estimation is designed for rotorcraft UAVs with a body-fixed thrust direction and three-axis attitude control.  The vehicle is modeled as an under-actuated system on the tangent bundle of the six-dimensional Lie group of rigid body motions, $\SE$. The proposed ESO scheme is developed based on the HC-FFTSD, which is similar to the STA used in sliding mode designs, to obtain fast finite-time stability with higher tunability of the settling time compared to other FTS schemes.  The Lyapunov stability analysis presented in this article for the ESO scheme proves the finite-time stability and robustness of the ESO on $\SE$. A set of numerical simulations are conducted.  The numerical simulation results present the stable performance of the FFTS-ESO scheme in estimating external force and torque disturbances acting on the UAV in different scenarios. The behavior of the FFTS-ESO is compared with two state-of-the-art observers for disturbance estimation. Using a realistic set of data for several simulated flight scenarios of a rotorcraft UAV, numerical simulations show that the FFTS-ESO, unlike the LESO and FxTSDO, is always stable and its convergence is robust to measurement noise and pose singularities. The proposed FFTS-ESO is implemented on the FCU of a multi-rotor UAV, with disturbance rejection control using feedback of disturbance estimates from the FFTS-ESO. The results validate the proposed FFTS-ESO experimentally and show the supremacy of the disturbance rejection control scheme over the original control scheme.

\appendix
\section{Proof of Lemma 5}
\begin{Proof*}
	{\em
		\textup{Represent $x$ as a linear combination of $\mu$ and $\nu$: }
		\begin{align}\label{eqn:x representation}
			x = c_1 \mu +c_2\nu,  
		\end{align}
		\textup{where $\nu$ is a vector perpendicular to $\mu$, such that $\mu\T\nu=0$. Next, define two non-zero scalars, $c_1, c_2$. Using \eqref{eqn:x representation}, express $Y$ in Lemma \ref{lem:Inequality Noise Robustness} in coordinates $(c_1,c_2)$: }
		\begin{align*}
			Y = \frac{c_1\mu+c_2\nu}{\left(c_1^2\|\mu\|^2 +c_2^2\|\nu\|^2\right)^\alpha} - \frac{(1+c_1)\mu+c_2\nu}{\left[(1+c_1)^2\|\mu\|^2 +c_2^2\|\nu\|^2\right]^\alpha}.
		\end{align*}
		\textup{Thereafter, we obtain its partial derivatives with respect to these coordinates:}
		\begin{align}\label{eqn:Y partial c1}
			\begin{split}
				\frac{\partial Y}{\partial c_1} &= \frac{\mu}{\left(c_1^2\|\mu\|^2 +c_2^2\|\nu\|^2\right)^\alpha}- \frac{2\alpha c_1\|\mu\|^2(c_1\mu+c_2\nu)}{\left(c_1^2\|\mu\|^2 +c_2^2\|\nu\|^2\right)^{\alpha+1}}  \\
				&- \frac{\mu}{\left[(1+c_1)^2\|\mu\|^2 +c_2^2\|\nu\|^2\right]^\alpha} + \frac{2\alpha( 1+c_1)\|\mu\|^2\left[(1+c_1)\mu+c_2\nu\right]}{\left[(1+c_1)^2\|\mu\|^2 +c_2^2\|\nu\|^2\right]^{\alpha+1}},
			\end{split}
		\end{align}
		\begin{align}\label{eqn:Y partial c2}
			\begin{split}
				\frac{\partial Y}{\partial c_2} &= \frac{\nu}{\left(c_1^2\|\mu\|^2 +c_2^2\|\nu\|^2\right)^\alpha} - \frac{2\alpha c_2\|\nu\|^2(c_1\mu+c_2\nu)}{\left(c_1^2\|\mu\|^2 +c_2^2\|\nu\|^2\right)^{\alpha+1}} \\
				&- \frac{\nu}{\left[(1+c_1)^2\|\mu\|^2 +c_2^2\|\nu\|^2\right]^\alpha} + \frac{2\alpha c_2\|\nu\|^2\left[(1+c_1)\mu+c_2\nu\right]}{\left[(1+c_1)^2\|\mu\|^2 +c_2^2\|\nu\|^2\right]^{\alpha+1}}.
			\end{split}
		\end{align}
		\textup{Thereafter, we employ the fact that the local maxima of $Y\Tp Y$ satisfy:}
		\begin{align*}
			\begin{split}
				\frac{\partial }{\partial c_1} (Y\Tp Y) = \frac{\partial }{\partial c_2} (Y\Tp Y) = 0, 
			\end{split}
		\end{align*}
		we obtain the following equivalent conditions for the maxima:
		\begin{align}
			&\nu\Tp \frac{\partial Y}{\partial c_1} = \mu\Tp \frac{\partial Y}{\partial c_2} = 0, 	\label{eqn:Y partial equation 01}\\
			&\mu\Tp \frac{\partial Y}{\partial c_1} = 0, 											\label{eqn:Y partial equation 02}\\
			&\nu\Tp \frac{\partial Y}{\partial c_2} = 0. 											\label{eqn:Y partial equation 03}
		\end{align}
		Substituting \eqref{eqn:Y partial c1} and \eqref{eqn:Y partial c2} into \eqref{eqn:Y partial equation 01}, we obtain:
		\begin{align}\label{eqn:Y partial equation 1a} 
			\begin{split}
				&\nu\Tp \frac{\partial Y}{\partial c_1} = \mu\Tp \frac{\partial Y}{\partial c_2} =  0, \\
				&\Longleftrightarrow -\frac{2\alpha c_1 c_2 \|\mu\|^2\|\nu\|^2 }{\left(c_1^2\|\mu\|^2 +c_2^2\|\nu\|^2\right)^{\alpha+1}} + \frac{2\alpha (1+c_1) c_2 \|\mu\|^2\|\nu\|^2}{\left[(1+c_1)^2\|\mu\|^2 +c_2^2\|\nu\|^2\right]^{\alpha+1}} =0, \\
				&\Longrightarrow  \ c_1\left[(1+c_1)^2\|\mu\|^2 +c_2^2 \|\nu\|^2\right]^{\alpha+1} = (1+c_1)\left[c_1^2\|\mu\|^2 +c_2^2 \|\nu\|^2\right]^{\alpha+1}, \\
			\end{split}
		\end{align}
		Substituting \eqref{eqn:Y partial c1} and \eqref{eqn:Y partial c2} into \eqref{eqn:Y partial equation 02}, we obtain:
		\begin{align}\label{eqn:Y partial equation 2}
			\begin{split}
				&\mu\Tp \frac{\partial Y}{\partial c_1} = 0, \\
				&\Longrightarrow  \frac{(1-2\alpha \|\mu\|^2 c_1^2)\|\mu\|^2}{\left(c_1^2\|\mu\|^2 +c_2^2\|\nu\|^2\right)^{\alpha+1}} - \frac{\left[ 1-2\alpha (1+c_1)^2\|\mu\|^2 \right]\|\mu\|^2}{\left[(1+c_1)^2\|\mu\|^2 +c_2^2\|\nu\|^2\right]^{\alpha+1}} = 0, \\
				&\Longleftrightarrow  (1+c_1)^2 = c_1^2, \Longleftrightarrow  c_1 = -\frac{1}{2}.
			\end{split}
		\end{align}
		Substituting \eqref{eqn:Y partial c1} and \eqref{eqn:Y partial c2} into \eqref{eqn:Y partial equation 03}, we obtain:
		\begin{align}\label{eqn:Y partial equation 3}
			\begin{split}
				&\nu\Tp \frac{\partial Y}{\partial c_2} = 0, \\
				&\Longrightarrow  \frac{(1-2\alpha \|\nu\|^2 c_2^2)\|\nu\|^2}{\left(c_1^2\|\mu\|^2 +c_2^2\|\nu\|^2\right)^{\alpha+1}} - \frac{(1-2\alpha \|\nu\|^2 c_2^2)\|\nu\|^2}{\left[(1+c_1)^2\|\mu\|^2 +c_2^2\|\nu\|^2\right]^{\alpha+1}} = 0, \\
				\Longleftrightarrow & (1+c_1)^2 = c_1^2, \Longleftrightarrow  c_1 = -\frac{1}{2}. \\
			\end{split}
		\end{align}
		\eqref{eqn:Y partial equation 1a} does not give a real solution for $\alpha \in ]0,1/2[  $. Thus, we conclude that the only solution to \eqref{eqn:Y partial equation 01}, \eqref{eqn:Y partial equation 02}, \eqref{eqn:Y partial equation 03} is given by $c_1 = -1/2, c_2 = 0$. 
		Thus, the only critical value of $Y\Tp Y$ is obtained when $x=-\mu/2$.
		Further, we conclude that the global maximum of $Y\Tp Y$ is at $x=-\mu/2$ because it is positive definite in $Y$. Therefore, we do not need an analysis of the Hessian matrix of $Y\T Y$ as a function of $(c_1,c_2)$. 
	} 
	\qed
\end{Proof*}

\section*{Acknowledgement}
The authors acknowledge support from the National Science Foundation award 2132799 and WindShape Corp.

\bibliographystyle{apalike} 
\bibliography{ref_CEP.bib}

\end{document}